\documentclass[11pt,twoside]{article} 
\usepackage{latexsym} 
\usepackage{enumitem}
\input{amssym.def} 
\input{amssym} 
\newfont{\fra}{eufm10 scaled 1095} 
\newfont{\Bb}{msbm10 scaled 1095} 
\newfont{\Bbg}{msbm10 scaled 1680} 
 
\newcommand\RR{{\mbox{\Bb R}}} 
 
\newcommand\ZZ{{\mbox{\Bb Z}}}

\newcommand\fb{{\frak b}} 
\newcommand\fg{{\frak{g}}} 
\newcommand\fh{{\frak h}} 
\newcommand\fri{{\frak i}} 
 
\newcommand\fl{{\frak l}} 
\newcommand\fm{{\frak m}} 
\newcommand\fn{{\frak n}}

\newcommand\fa{{\frak a}} 
\newcommand\fd{{\frak d}} 
\newcommand\fr{{\frak r}}

\newcommand\fz{{\frak z}} 

\newcommand\cH{{\cal H}}

\newcommand\cN{{\cal N}}
\newcommand\ph{\varphi} 
 
\newcommand\eps{\varepsilon} 

\newcommand{\fsl}{\mathop{{\frak s \frak l}}} 
 
\newcommand\fspin{\frak s \frak p \frak i \frak n}
 
\newcommand{\fso}{\mathop{{\frak s \frak o}}}

\newcommand{\End}{\mathop{{\rm End}}} 
\newcommand{\Mat}{\mathop{{\rm Mat}}}

\newcommand{\GL}{\mathop{{\rm GL}}} 
\newcommand{\SO}{\mathop{{\rm SO}}}
\newcommand{\SU}{\mathop{{\rm SU}}}
\newcommand{\Sp}{\mathop{{\rm Sp}}}
\newcommand{\grO}{{{\rm O}}}
\newcommand{\Spin}{\mathop{{\rm Spin}}} 
\newcommand{\G}{\mathop{{{\rm G}_{2(2)}}}} 
 
\newcommand{\Aut}{\mathop{{\rm Aut}}} 

\newcommand{\Id}{{{\rm id}}} 
\newcommand{\ad}{{{\rm ad}}} 
\newcommand{\tad}{{\widetilde{\rm ad}}} 
\newcommand{\tr}{\mathop{{\rm tr}}}

\newcommand{\Ric}{\mathop{{\rm Ric}}}

\newcommand{\Cl}{{\cal{C}}}

\newcommand{\Span}{{{\rm span}}} 
\newcommand{\mod}{\mathop{{\rm mod}}} 
\newcommand{\proj}{{{\rm pr}}} 
\newcommand\ip{{\langle\cdot \,,\cdot \rangle}} 
\newcommand\lb{{[\cdot\,,\cdot]}}

\newcommand\dd{\fd_{\alpha,\gamma}(\fl,\theta_\fl,\fa)} 
\newcommand\ddm{\fd_{\alpha,\gamma}(\fl,\theta_\fl,b_\fm,\fa)} 
\newcommand\HQ{{{\cal H}_Q^2(\fl,\fa)}}
\newcommand\HQb{{{\cal H}_Q^2(\fl,\theta_\fl,b_\fm,\fa)}} 
\newcommand\ZQb{{{\cal Z}_Q^2(\fl,\theta_\fl,b_\fm,\fa)}} 
\newcommand\CQtheta{{{\cal C}_Q^1(\fl,\fa)_+}} 
\newcommand\HQtheta{{{\cal H}_Q^2(\fl,\theta_\fl,\fa)}} 
\newcommand\ZQtheta{{{\cal Z}_Q^2(\fl,\fa)_+}} 

\newcommand\proof{{\sl Proof. }} 
\newcommand{\qed}{\hspace*{\fill}\hbox{$\Box$}\vspace{2ex}} 
\newcommand{\benur}{\begin{enumerate}[label=(\roman*)]}
\newcommand{\la}{\langle}
\newcommand{\ra}{\rangle}
\newtheorem{theo}{Theorem}[section] 
\newtheorem{pr}[theo]{Proposition}
\newtheorem{de}[theo]{Definition}

\newtheorem{re}[theo]{Remark}
\newtheorem{co}[theo]{Corollary}
\newtheorem{lm}[theo]{Lemma}
\begin{document} 
\title{Indefinite symmetric spaces with $\G$-structure} 
\author{Ines Kath}
\maketitle 
\begin{abstract} \noindent We determine all indecomposable pseudo-Riemannian symmetric spaces of signature $(4,3)$ whose holonomy is 
contained in $\G\subset SO(4, 3)$.
\end{abstract}
\section{Introduction}
The compact Lie group ${\rm G}_2$ lies on the list of holonomy groups of irreducible Riemannian manifolds. Each Riemannian manifold whose holonomy group is contained in ${\rm G}_2$ is Ricci-flat. In particular, if the holonomy group of a Riemannian symmetric space $M$ is contained in ${\rm G}_2$, then $M$ must be flat. In contrast, we will see that there exist indecomposable indefinite symmetric spaces of signature $(4,3)$ whose holonomy is contained in the split real form $\G\subset \SO(4,3)$. The same effect appears also with other holonomy groups known from Riemannian geometry, e.g, there are Ricci-flat hermitian symmetric spaces (see \cite{KO2} for index 2) and hyper-K\"ahler symmetric spaces, see~\cite{AC1, KO3}. Of course, in all cases the holonomy group of the symmetric space is properly contained in $\G$, $\SU(p,q)$, and $\Sp(p,q)$, respectively. 

Let us recall the definition of the group $\G$. We consider the generic 3-form 
\begin{equation}\label{Eomega}
\omega_0=\sqrt2 (\sigma^{127}+\sigma^{356})-\sigma^4\wedge(\sigma^{15}+\sigma^{26}-\sigma^{37})\,.
\end{equation}
on $\RR^7$, where $\sigma^1,\dots,\sigma^7$ denotes the dual basis of the standard basis and $\sigma^{ij}:=\sigma^i\wedge\sigma^j$, $\sigma^{ijk}:=\sigma^i\wedge\sigma^j\wedge \sigma^k$. The group $\GL(7)$ acts on the space of 3-forms on $\RR^7$ and we define $\G\subset \GL(7)$ to be the stabiliser of $\omega_0$.  Then $\G$ is a non-compact group of dimension 14. It is contained in the orthogonal group with respect to the scalar product  
$2\sigma^1\cdot\sigma^5+2\sigma^2\cdot \sigma^6+2\sigma^3\cdot \sigma^7- (\sigma^4)^2$, which has signature $(4,3)$. It is connected and its fundamental group is $\ZZ_2$.
There are other nice characterisations of this group, e.g., $\G$ is the stabiliser of a non-isotropic element of the real spinor representation of $\Spin(4,3)$ and it can also be understood as the stabiliser of a cross product on $\RR^{4,3}$.

For a Lie group $G\subset \grO(p,q)$, a $G$-structure on a pseudo-Riemannian manifold $(M,g)$ of signature $(p,q)$ is a reduction of the bundle $P_{SO(p,q)}$ of orthonormal frames to a $G$-bundle $P_G\subset P_{SO(p,q)}$. It is called parallel if $P_G$ is parallel with respect to the Levi-Civita connection. The existence of a parallel $G$-structure is equivalent to the fact that the holonomy group of $(M,g)$ is contained in $G$. Pseudo-Riemannian symmetric spaces with parallel $G$-structures are studied, e.g., for $G=U(r,s)\subset\SO(2r,2s)$ (pseudo-Hermitian symmetric spaces) in \cite{KO2, KOesi}, for $G=\Sp(r,s)\subset\SO(4r,4s)$ (hyper-K\"ahler symmetric spaces) in \cite{AC1,KO3} and for $G=\Sp(r,s)\cdot\Sp(1)\subset\SO(4r,4s)$ (quaternionic K"ahler case) in \cite{AC2}. Here we turn to $G=\G\subset \SO(4,3)$. According to the above remarks on the definition of $\G$ a parallel $\G$-structure can equivalently be defined as a parallel 3-form on $M$ that equals $\omega_0$ with respect to a suitable local frame. Analogously, parallel non-isotropic spinor fields or parallel cross-products can be used to define parallel $\G$-structures. 

The aim of the paper is to classify all parallel $\G$-structures on (simply-connected) indecomposable pseudo-Riemannian symmetric spaces of signature $(4,3)$.  Since any symmetric space $(M,g)$ that admits a $\G$-structure is Ricci-flat its transvection group is solvable. In particular, its holonomy group is a proper subgroup of $\G$ as already remarked above.
Since the transvection group is solvable, the symmetric triple $(\fg,\theta,\ip)$ that is associated with $(M,g)$ is  the quadratic extension of a Lie algebra with involution $(\fl,\theta_\fl)$ by an orthogonal $(\fl,\theta_\fl)$-module $\fa$ as described in \cite{KO2}. Hence we can apply the structure theory for such extensions developed in \cite{KO2}.

Henceforth we will use the following notation. If $\fg$ is a Lie algebra and $\theta$ is an involution on $\fg$, then $\fg_\pm$ will denote the eigenspace of $\theta$ with eigenvalue $\pm1$. 

Now let $(M,g)$ be a symmetric space with a parallel $\G$-structure and $(\fg,\theta,\ip)$ the symmetric triple associated with $(M,g)$. We think of the $\G$-structure as a parallel cross product $b$ on $(M,g)$. Then $b$ corresponds to a $\fg_+$-invariant cross product on $\fg_-$, which we will also denote by $b$. Up to coverings, $(M,g)$ and its $\G$-structure are determined by the symmetric triple $(\fg,\theta,\ip)$ and $b:\fg_-\times\fg_-\rightarrow \fg_-$. As remarked above $(\fg,\theta,\ip)$ has the structure of a quadratic extension. Moreover, this structure is uniquely determined by the so-called canonical isotropic ideal $\fri\subset\fg$ as explained in \cite{KO2}. Our first step will be to prove that $\fri_-=\fri\cap\fg_-$ is three-dimensional provided the symmetric triple $(\fg,\theta,\ip)$ is indecomposable.  Then we show that $\fri_-$ is invariant under $b$. Moreover, after a suitable choice of a section $\fl_-:=\fg_-/\fri_-^\perp\hookrightarrow \fg_-$, $b$ also restricts to a bilinear map on $\fl_-$. This will give an additional structure on $\fl$, which will allow to determine first $\fl$ and than the second quadratic cohomology of $\fl$ with coefficients in $\fa:=\fri^\perp/\fri$. This cohomology set is in bijection with the equivalence classes of quadratic extensions of $(\fl,\theta_\fl)$ by $\fa$, which will yield a classification. The final result is  a list of all indecomposable symmetric triples with $\G$-structure, see Section~\ref{Sfinal}.  It is formulated in a self-contained manner, i.e., without requiring the knowledge of the results in \cite{KO2}.

I would like to thank Martin Olbrich for his interest in this project and for his valuable comments.

\section{The spinor representation in signature $(4,3)$}
Let $V$  be an oriented 7-dimensional vector space endowed with a scalar product $\ip$ of signature $(4,3)$. 
Let $b_1,\dots,b_7$ be a basis of $V$ and denote the dual basis of $V^*$ by $\sigma^1,\dots,\sigma^7$. We will say that $b_1,\dots b_7$ is a Witt basis if the scalar product on $V$ equals 
$$2\sigma^1\cdot\sigma^5+2\sigma^2\cdot \sigma^6+2\sigma^3\cdot \sigma^7- (\sigma^4)^2.$$
The standard example is the vector space $\RR^7$ together with the the scalar product $\ip:=-(dx_4)^2+2\sum_{i=1}^3 dx_i dx_{i+4}$, which we will denote by $\RR^{4,3}$. The standard basis $e_1,\dots,e_7$ of $\RR^7$ is a Witt basis of $\RR^{4,3}$. 
     
The Clifford algebra $\Cl_V$ is the unital associative algebra that is generated by the elements of $V$ subject to the relations $uv+vu=-2\la u,v\ra$. It is isomorphic to a sum of real matrix algebras, more exactly, $\Cl_V\cong \Mat(8,\RR)\oplus\Mat(8,\RR)$. This gives us two inequivalent irreducible representations of $\Cl_V$ on $\RR^8$. Let $\hat b_1,\dots,\hat b_7$ be a positively oriented orthonormal basis of $V$ and denote by $\alpha_0:=\hat b_1\cdot \dots\cdot \hat b_7$ the volume element in $\Cl_V$. Then the equivalence classes of irreducible representations of $\Cl_V$ differ by the action of $\alpha_0$, which is either the identity or minus identity.  If $\alpha_0$ acts by $\Id$ we will say that the representation is of type one, otherwise it is of type two. Let $\Delta_V$ denote the equivalence class of representations of type one.

Now consider the group 
$$\Spin(V):= \la uv\mid u,v\in V,\ \la u,u\ra=\pm1, \la v,v\ra=\pm1\ra  \subset {\Cl}_V.$$
If we restrict the two inequivalent representations of $\Cl_V$ to $\Spin(V)$ we obtain two equivalent irreducible representations of $\Spin(V)$ on $\RR^8$, which are called spinor representation of $\Spin(V)$. 
The Lie algebra $\fspin(V)$ of $\Spin(V)$ equals 
$$\fspin(V)=\Span\{ \hat b_i\hat b_j\mid 1\le i < j \le 7\}\subset\Cl_V.$$

There is a two-fold covering map $\lambda:\Spin(V)\rightarrow \SO(V)$, which is defined by $\lambda(a)(v)=ava^{-1}$ for $a\in\Spin(V)$ and $v\in V$. In particular, $\lambda$ induces an isomorphism $\lambda_*: \fspin(V)\rightarrow \fso(V)$. The inverse of this isomorphism is given by
\begin{equation}\label{Etilde}
\fso(V)\ni A\longmapsto \ \tilde A=\frac14 \sum_{i=1}^3\Big(b_i A(b_{i+4}) + b_{i+4}A(b_i)\Big)- \frac14 b_4A(b_4)\in\fspin(V).
\end{equation}

There exists an inner product $\ip_\Delta$ of signature $(4,4)$ on $\Delta_V$ satisfying
\begin{equation} \label{Eipd}
\la X \cdot \ph,\psi\ra_\Delta+\la \ph,X\cdot \psi\ra_\Delta=0
\end{equation}
for all $X\in V$.
 It is uniquely determined up to multiplication with a real number $\lambda\not=0$. Independently of the choice of $\ip_\Delta$ we can speak of isotropic and non-isotropic spinors and of pairs of orthogonal spinors.

Let us give explicit formulas in the case, where $V=\RR^{4,3}$. We denote the Clifford algebra of $\RR^{4,3}$ by $\Cl_{4,3}$  and the spin group by $\Spin(4,3)$. Let $s_1,\dots,s_8$ denote the standard basis of $\RR^8$. We define an algebra homomorphism $\Phi:\Cl_{4,3}\to \Mat(8,\RR)$ by 
\begin{eqnarray}\label{EPhi}
(1/\sqrt2)\cdot\Phi(e_1)&:& s_1\mapsto s_8,\ s_2\mapsto s_7,\ s_3\mapsto -s_6,\ s_4\mapsto -s_5\nonumber\\
(1/\sqrt2)\cdot\Phi(e_2)&:& s_1\mapsto -s_3,\ s_2\mapsto s_4,\ s_7\mapsto s_5,\ s_8\mapsto -s_6 \nonumber \\
(1/\sqrt2)\cdot\Phi(e_3)&:& s_2\mapsto -s_1,\ s_4\mapsto -s_3,\ s_5\mapsto s_6,\ s_7\mapsto s_8\nonumber \\
(1/\sqrt2)\cdot\Phi(e_5)&:& s_5\mapsto s_4,\ s_6\mapsto s_3,\ s_7\mapsto -s_2,\ s_8\mapsto -s_1\\
(1/\sqrt2)\cdot\Phi(e_6)&:& s_3\mapsto s_1,\ s_4\mapsto -s_2,\ s_5\mapsto -s_7,\ s_6\mapsto s_8 \nonumber\\
(1/\sqrt2)\cdot\Phi(e_7)&:& s_1\mapsto s_2,\ s_3\mapsto s_4,\ s_6\mapsto -s_5,\ s_8\mapsto -s_7 \nonumber \\
\Phi(e_4)&:& s_i\mapsto s_i,\ i=1,4,6,7,\quad s_j\mapsto -s_j,\ j=2,3,5,8, \nonumber
\end{eqnarray}
where all basis vectors that are not mentioned are mapped to zero.
This will give us a  representation of type one of $\Cl_{4,3}$, which we will denote by $\Delta_{4,3}$.
The inner product $$\ip_\Delta=2\sum_{i=1}^4 dx_i dx_{i+4}$$ on $\Delta_{4,3}\cong \RR^8$ satisfies (\ref{Eipd}). 
\begin{pr} \label{Ciso2}
Let $c\in\RR$ be fixed.
The group $\Spin(V)$ acts transitively on 
$$\Delta(c):=\{ \psi\in \Delta_{V}\mid \langle \psi,\psi\rangle_\Delta =\pm c\}.$$ 
If $c\not=0$, then the stabiliser of an element $\psi\in \Delta(c)$ is isomorphic to $\G$. 

Furthermore, $\Spin(V)$ acts transitively on 
$$\Delta(c,0):=\{(\psi,\ph)\in \Delta_{V}\times \Delta_{V}\mid \ph\perp \psi, \ \la \psi,\psi\ra_\Delta=\pm c,\  \langle \ph,\ph\ra_\Delta=0,\ \ph\not=0\}.$$
The Lie algebra $\fh(\psi,\ph)$ of the stabiliser of an element $(\psi,\ph)\in\Delta(c,0)$ is isomorphic to $\fsl(2,\RR)\ltimes \fm$, where $\fm$ is a central extension of the 3-dimensional Heisenberg algebra by a two-dimensional vector space. If $\fa\subset \fh(\psi,\ph)$ is an abelian subalgebra, then $\dim\fa\le 3$.
\end{pr}
\proof The fact that $\Spin(V)$ acts transitively on $\Delta(c)$ with stabiliser isomorphic to $\G$ is well-known, for a proof see, e.g., \cite{KG2}. Furthermore, in \cite{KG2}, Prop. 2.3. it is shown that the unity component $\Spin^+(4,3)$ of $\Spin(4,3)$ acts transitively on the Stiefel manifolds
$$V(\eps_1,\eps_2,\eps_3):=\{(\psi_1,\psi_2,\psi_3)\in\Delta_{4,3}\mid \la \psi_i,\psi_j\ra_\Delta=\eps_i\delta_{ij}\}$$ 
for $\eps_i=\pm1$, $i=1,2,3$. This shows that $\Spin(4,3)$ acts transitively on $\Delta(c,0)$ since, for a given $(\psi,\ph)\in\Delta(c,0)$ we can write $\ph=\psi_1+\psi_2$ such that $\psi,\psi_1,\psi_2$ are orthogonal and $\la\psi_1,\psi_1\ra_\Delta=-\la\psi_2,\psi_2\ra_\Delta=c$. 

The Lie algebra $\fh:=\fh(s_1+s_5,s_6)\subset\fspin(4,3)$   is spanned by
\begin{eqnarray*}
&e_1e_5-e_2e_6,\ e_1e_6,\ e_2e_5,\ Z_1:=e_1e_3,\ Z_2:=e_2e_3,&\\ &N_1:=e_3e_5-\sqrt2 e_2e_4,\ N_2:=e_3e_6+\sqrt2 e_1e_4,\ N_3:=e_1e_2+\sqrt2 e_3e_4. &
\end{eqnarray*}
Let us first analyse the structure of the Lie algebra $\fh$. Note that $e_1e_6,\ e_2e_5$ and  $e_1e_5-e_2e_6$ span a subalgebra of $\fh$ that is isomorphic to $\fsl(2,\RR)$. Moreover, $Z_1,\ Z_2,\ N_1,\ N_2,\ N_3$  span an ideal $\fm$ of $\fh$, i.e., 
\begin{equation}\label{Esl}
\fh=\fsl(2,\RR)\ltimes \fm.
\end{equation}
The centre of $\fm$ equals $\fz(\fm):=\Span\{Z_1,Z_2\}$. Hence, as a vector space, $\fm$ is the direct sum of $\fz(\fm)$ and  $\fn:=\Span\{N_1,N_2,N_3\}$. Since
\begin{equation}\label{EN}
[N_1,N_2]=-4N_3 ,\ [N_1,N_3]=6Z_2  ,\ [N_2,N_3]= -6Z_1
\end{equation}
$\fm$ is a central extension of the three-dimensional Heisenberg algebra $\fh(1)$ by $\fz(\fm)$. The representation of $\fsl(2,\RR)\subset \fh$ on $\fm$ decomposes into a one-dimensional trivial representation spanned by $N_3$ and two standard representations spanned by $Z_1,Z_2$ and $N_1,N_2$, respectively. Now let $\fa\subset\fh$ be an abelian subalgebra. 

Let us first consider the case where $\fa\subset\fm$. Let $\proj_\fn$ denote the projection from $\fm=\fz(\fm)\oplus\fn$ to $\fn$. Then $\proj_\fn(\fa)$ is also an abelian subalgebra of $\fm$. Using (\ref{EN}) we see that $\proj_\fn(\fa)$ is at most one-dimensional. Hence $\fa$ is at most three-dimensional and $\dim\fa=3$ holds if and only if $\fa=\fz(\fm)\oplus \RR\cdot N$ for some $N\not=0$ in $\fn$.

Now suppose $\fa\not\subset\fm$. Since $\fsl(2,\RR)$ does not contain an abelian subalgebra of dimension two the projection of $\fa$ to $\fsl(2,\RR)$ with respect to the decomposition (\ref{Esl}) is one-dimensional. Thus $\fa=\RR(B+M)\oplus\fa_\fm$ for some $B\in \fsl(2,\RR)$ and $M\in \fm$, where $\fa_\fm$ is an abelian subalgebra of $\fm$. We have to show that $\dim \fa_\fm\le2$ holds. Assume $\dim\fa_\fm=3$. Then the above considerations show $\fa_\fm=\fz(\fm)\oplus \RR\cdot N$ for some $N\not=0$ in $\fn$. Since $\fa$ is abelian $B\not=0$ has to act trivially on $\fz(\fm)\subset\fa_\fm$, which is a contradiction. \qed 
\begin{pr}\label{Pmult}
 \begin{enumerate} 
\item If $\psi\in \Delta_{V}$ is non-isotropic, then $V\ni X\mapsto X\cdot \psi\in \psi^\perp$ is an isomorphism.
\item The map
$$\Delta_V\ni \ph\longmapsto U_\ph:=\{ X\in V\mid X\cdot \ph=0\}\subset V$$
induces a bijection from the set 
$\{\ph \in \Delta_V\mid \la\ph,\ph\ra=0,\ \ph\not=0\}/\RR$ of projective isotropic spinors to the set of  3-dimensional isotropic subspaces of $V$. If $\ph\in\Delta_V$ is isotropic, then $U_\ph^\perp\cdot \ph\subset\RR\cdot \ph$.
\end{enumerate}
\end{pr}
\proof By Prop.~\ref{Ciso2} we may assume $\psi=s_1+s_5\in\Delta_{4,3}$ and $\ph=s_6\in\Delta_{4,3}$ and the assertion follows from Equation~(\ref{EPhi}).\qed
\begin{de} For the time being, consider $V$ without orientation. A 3-form $\omega$ on $V$ is called nice if there is a Witt basis $b_1,\dots,b_7$ of $V$ such that  
\begin{equation}\label{Enochmal}
\omega=\sqrt2 (\sigma^{127}+\sigma^{356})-\sigma^4\wedge(\sigma^{15}+\sigma^{26}-\sigma^{37})\,.
\end{equation}
with respect to the dual basis $\sigma^1,\dots,\sigma^7$. 
\end{de}
The stabiliser of a nice 3-form $\omega$ is isomorphic to $\G$. In particular, $\omega$ induces an orientation on $V$, since $\G$ is connected.
\begin{de}\label{Db}
Consider again $V$ without orientation. A bilinear map $b:V\times V\rightarrow V$ is called a cross product if
\benur
\item $b(X,Y)=-b(Y,X)$,
\item $\langle X,b(X,Y)\rangle=0$,
\item $b(X,b(X,Y))=-\langle X,X\rangle Y +\langle X,Y\rangle X$.
\end{enumerate}
\end{de}
The following proposition is proven in \cite{KG2}.
\begin{pr}\label{Pobp} \begin{enumerate}
\item The map 
$${\cal B}\longrightarrow \textstyle{\bigwedge}_\sharp^3(V^*), \quad b\longmapsto \omega_b:=\langle\, \cdot\,, b(\cdot,\cdot)\ra$$
is a bijection between the set ${\cal B}$ of cross products and the set $\bigwedge_\sharp^3(V^*)$ of nice 3-forms on $V$.

\item Now let $V$ be oriented and denote by ${\cal B}^+$ the set of cross products $b$ for which the orientation induced by $\omega_b$ coincides with the orientation of $V$. Then the map 
$$\Delta_*:=\{\psi\in\Delta_V\mid \la\psi,\psi\ra_\Delta\not=0\}\longrightarrow {\cal B}^+,\quad \psi\longmapsto b_\psi$$
defined by
$$ XY\cdot \psi + \langle X,Y\rangle \psi = b_\psi(X, Y)\cdot\psi$$
induces a bijection from the set $P(\Delta_*):=\Delta_*/\RR$ of projective non-isotropic spinors in $\Delta_V$ to ${\cal B}^+$.
\end{enumerate}
\end{pr}
If we consider, in particular, $\psi=s_1+s_5\in\Delta_{4,3}$, then the 3-form that is associated with $\psi$ according to Proposition~\ref{Pobp} equals $\omega_0$ as defined in~(\ref{Eomega}).
\section{$\G$-structures on symmetric spaces}
\subsection{Symmetric spaces and symmetric triples}
 Before we start let us introduce the following convention. If $\fg$ is a Lie algebra and $\theta$ is an involutive automorphism on $\fg$, then we denote the eigenspaces of $\theta$ with eigenvalues $1$ and $-1$ by $\fg_+$ and $\fg_-$, respectively. 
 
Let $M$ be a (pseudo-Riemannian) symmetric space and choose a base point $u\in M$. Then $M=G/G_+$, where $G$ is the transvection group of $M$ and $G_+\subset G$ is the stabiliser of $u\in M$. The conjugation by the reflection of $M$ at $u$ induces an involution on $G$ and therefore also on $\fg$. We denote this involution by $\theta$. If $\fg=\fg_+\oplus\fg_-$ is the decomposition of $\fg$ into eigenspaces of $\theta$, then $\fg_+$ is the Lie algebra of $G_+$ and the vector space $\fg_-$ can be identified with $T_uM$. Moreover, $G_+$ equals the holonomy group of $(M,g)$ and its adjoint representation on $\fg_-$ is the holonomy representation. The scalar product on $\fg_-\cong T_uM$ has a unique extension to an $\ad(\fg)$-invariant non-degenerate inner product $\ip$ on $\fg$ such that $\fg_+\perp\fg_-$. In particular, the triple $(\fg, \theta,\ip)$ consists of a metric Lie algebra $(\fg,\ip)$ and an isometric involution $\theta:\fg\rightarrow\fg$. Moreover, $[\fg_-,\fg_-]=\fg_+$. Any triple $(\fg, \theta,\ip)$ satisfying these properties will be called a symmetric triple. The above described assignment of a symmetric triple to a symmetric space gives a bijection from the set of simply-connected symmetric spaces to the set of symmetric triples. Isometry classes of simply-connected symmetric spaces correspond to isomorphism classes of symmetric triples. Furthermore, a symmetric space is indecomposable if and only if the associated symmetric triple is indecomposable, i.e., if it is not a non-trivial direct sum of two symmetric triples.

The signature of a symmetric triple $(\fg, \theta,\ip)$ is defined as the signature of $\ip$ restricted to $\fg_-$, which equals the signature of $M$.
\subsection{$\G$-structures on symmetric spaces and symmetric triples}
\begin{de} A $\G$-structure on a pseudo-Riemannian manifold $(M,g)$ of signature $(4,3)$ is a section $\omega\in\Omega^3(M)$ such that $\omega_x$ is a nice 3-form on $T_xM$ for each $x\in M$.
\end{de}
\begin{de}\label{Dg2}
A  $\G$-structure $\omega$ on a symmetric triple $(\fg,\theta,\ip)$  of signature $(4,3)$ is a nice $\fg_+$-invariant 3-form $\omega$ on $\fg_-$. 
\end{de}
According to Prop.~\ref{Pobp}, it can equivalently be considered as a  $\fg_+$-invariant cross product $b$ on $\fg_-$ or as a pair $({\cal O},[\psi])$, where ${\cal O}$ is an orientation on $\fg_-$, $\psi\in \Delta_{\fg_-}$ is a $\fg_+$-invariant non-isotropic element of the representation $\Delta_{\fg_-}$ of $\Cl(\fg_-)$ of type one and $[\psi]$ is the projective spinor represented by $\psi$. 
\begin{de}
Two symmetric triples with $\G$-structure $(\fg_i,\theta_i,\omega_i,\ip_i)$, $i=1,2$, are called isomorphic if there is an isomorphism $\phi:(\fg_1,\theta_1,\ip_1)\rightarrow(\fg_2,\theta_2,\ip_2)$ of symmetric triples satisfying $\phi*\omega_2=\omega_1$.
\end{de}
The following proposition is a consequence of the holonomy principle.
\begin{pr} There is a 1-1-correspondence between parallel $\G$-structures on $(M,g)$ and $\G$-structures on the associated symmetric triple $(\fg,\theta,\ip)$.
\end{pr}
Since the stabiliser of a nice 3-form is isomorphic to $\G$ the existence of a parallel $\G$-structure implies that the holonomy group is contained in $\G$.
\begin{pr}\label{Pric}
If a symmetric space $M$ admits a parallel $\G$-structure, then its transvection group is solvable.
\end{pr}
\proof If a pseudo-Riemannian manifold of signature $(4,3)$ has  a parallel $\G$-structure, then it admits a spin structure and a parallel non-isotropic spinor field $\psi$. In particular, $Ric(X)\cdot\psi=0$ holds for any vector field $X$, see, e.g., \cite{BFGK}, hence $M$ is Ricci-flat. On the other hand, the Ricci tensor on $T_uM\cong\fg_-$ is given by the Killing form $\kappa_\fg$, more exactly, $\Ric(X,Y)=-(1/2)\cdot \kappa_\fg(X,Y)$. Hence, $\kappa_\fg(\fg_-,\fg_-)=0$. Moreover, $\kappa_\fg(\fg_+,\fg_+)=\kappa_\fg([\fg_-,\fg_-],\fg_+)\subset \kappa_\fg(\fg_-,\fg_-)=0$. Furthermore, $\kappa_\fg(\fg_+,\fg_-)=0$ since $\kappa_\fg$ is invariant under automorphisms, in particular, under $\theta$. Thus $\kappa_\fg=0$. Consequently, $\fg$ is solvable by Cartan's first criterion. \qed

In particular, this shows that the holonomy group of a symmetric space with $\G$-structure solvable, thus it is always properly contained in $\G$. Here we want to assume, that the holonomy group does not become `too small'. More exactly, we will consider only indecomposable symmetric spaces with $\G$-structure. 
\section{Quadratic extensions and $\G$-structures}
\subsection{Quadratic extensions}
Since the existence of a $\G$-structure on a symmetric triple $(\fg,\theta,\ip)$ implies that $\fg$ is solvable, the Lie algebra $\fg$  does not have simple ideals, thus the theory of quadratic extension applies. This theory is developed in \cite{KO2}. There it is proven that any symmetric triple $(\fg,\theta,\ip)$ without simple ideals has the structure of a canonically determined admissible quadratic extension of some proper Lie algebra with involution $(\fl,\theta_\fl)$ by an orthogonal $(\fl,\theta_\fl)$-module $\fa:=(\rho,\fa,\ip_\fa,\theta_\fa)$. Here the condition {\it proper} for $(\fl,\theta_\fl)$ means that $[\fl_-,\fl_-]=\fl_+$ holds. Moreover, any admissible quadratic extension of $(\fl,\theta_\fl)$ by $\fa$ is  equivalent to one of the following kind called standard model. It is denoted by $\dd$ and is obtained in the following way. As a symmetric triple it equals  $(\fd,\theta,\ip)$, where $\fd,\theta$ and $\ip$ are defined as follows. Assume that $(\alpha,\gamma)\in\ZQtheta$ is an admissible cocycle, i.e., $[\alpha,\gamma]\in\HQtheta_\sharp\subset\HQtheta$. Note that this assumption includes the condition that the representation of $\fl$ on $\fa$ is semisimple. Then $\fd$ equals the vector space $\fl^*\oplus\fa\oplus\fl$, the  inner product $\ip$ and an involutive endomorphism $\theta$ on $\fd$ are given by
\begin{eqnarray*}
 \langle Z+A+L,Z'+A'+L'\rangle&:=& \langle A,A'\rangle_\fa
+Z(L') +Z'(L) \\
\theta(Z+A+L)&:=& \theta_\fl ^*(Z)+\theta_\fa(A)+\theta_\fl(L)
\end{eqnarray*}
for $Z,\,Z'\in \fl^*$, $A,\,A'\in \fa$ and
$L,\,L'\in \fl$. Furthermore, the Lie bracket $\lb:\fd\times\fd\rightarrow \fd$ is defined by $[\fl^*,\fl^*\oplus\fa] =0$ and
\begin{eqnarray}
\ [L,L'] &=& \gamma(L,L',\cdot) +\alpha(L,L')+[L,L']_\fl\nonumber\\
\ [L,A] &=& \rho(L)(A) - \langle A,\alpha(L,\cdot)\rangle\nonumber\\
\ [L,Z]& = & \ad ^*(L)(Z)\label{Eduard}\\
\ [A,A']&=&\langle\rho(\cdot)(A),A'\rangle\nonumber
\end{eqnarray}
for $Z\in \fl^*$, $A,\,A'\in \fa$ and
$L,\,L'\in \fl$.
We identify the vector space $\fd/\fl^*$ with $\fa\oplus\fl$ and denote  by $i: \fa \rightarrow \fa \oplus \fl$ the injection and by $p:\fa\oplus \fl \rightarrow \fl$  the projection. Then  $(\fd,\fl^*,i,p)$ is an admissible quadratic extension of $(\fl,\theta_\fl)$ by $\fa$. We denote this quadratic extension as well as the underlying symmetric triple by $\dd$.

Recall that for any metric Lie algebra $\fg$ the canonical isotropic ideal is defined \cite{KO1}.
Since the considered quadratic extension $\dd$ is admissible the canonical isotropic ideal $\fri$  coincides with $\fl^*$.

For admissible cocycles $(\alpha_1,\gamma_1)$ and $(\alpha_2,\gamma_2)$ in $\ZQtheta$ the quadratic extensions $\fd_{\alpha_i,\gamma_i}(\fl,\theta_\fl,\fa)$, $i=1,2$, are equivalent if and only if $[\alpha_1,\gamma_1]=[\alpha_2,\gamma_2]\in\HQtheta_\sharp$.

Let $(\fl_i,\theta_{\fl_i})$, $i=1,2$, be Lie algebras with involution and let $\fa_i$ be orthogonal $(\fl_i,\theta_{\fl_i})$-modules. An isomorphism of triples $(S,U):(\fl_1,\theta_{\fl_1},\fa_1)\rightarrow(\fl_2,\theta_{\fl_2},\fa_2)$ consists of an isomorphism $S:(\fl_1,\theta_{\fl_1})\rightarrow (\fl_2,\theta_{\fl_2})$ of Lie algebras with involution and an isometry $U:\fa_2\rightarrow\fa_1$ such that $U\circ\rho_2 (S (L)) = \rho_1 (L) \circ U$ holds for all $L\in\fl_1$.

Let $(\fl,\theta_\fl)$ and $(\fl,\theta_{\fl'})$ be proper Lie algebras with involution. For $[\alpha,\gamma]\in\HQtheta_\sharp$ and $[\alpha',\gamma']\in\cH_Q^2(\fl',\theta_{\fl'},\fa')_\sharp$ the symmetric triples $\fd_{\alpha,\gamma}(\fl,\theta_\fl,\fa)$ and $\fd_{\alpha',\gamma'}(\fl',\theta_{\fl'},\fa')$ are isomorphic if and only if there is an isomorphism of triples $(S,U):(\fl,\theta_{\fl},\fa)\rightarrow(\fl',\theta_{\fl'},\fa')$ such that $(S,U)^*[\alpha',\gamma']=[\alpha,\gamma]$.
\begin{pr} \label{Pricl}
A symmetric space $(M,g)$ that is associated with a quadratic extension of $(\fl,\theta_\fl)$ by $\fa$ is Ricci-flat if and only if the trace form of the representation $\rho$ of $\fl$ on $\fa$ and the Killing form $\kappa_\fl$ of $\fl$ are related by $t_\rho=-2\kappa_\fl$.
\end{pr}
\proof In the proof of Prop.~\ref{Pric} we have already seen that $(M,g)$ is Ricci-flat if and only if $\kappa_\fg=0$ holds for the associated symmetric triple $(\fg,\theta,\ip)$. If $\fg$ is isomorphic to a quadratic extension of $(\fl,\theta_\fl)$ by $\fa$ then $\kappa_\fg(X_1,X_2)=0$ for $X_1,X_2\in\fl^*\oplus\fa$ and $\kappa_\fg(L_1,L_2)=2\kappa_\fl(L_1,L_2) +t_\rho(L_1,L_2)$ for all $L_1,L_2\in\fl$
by Equations (\ref{Eduard}).\qed
\subsection{The dimension of $ \fl_-$}\label{S32}
We already know that any symmetric triple that admits a $\G$-structure is a quadratic extension of a Lie algebra with involution $(\fl,\theta_\fl)$ by an orthogonal $(\fl,\theta_\fl)$-module $\fa$. Now let $(\fg,\theta,\ip)$ be a symmetric triple of signature $(4,3)$ that is isomorphic to the admissible quadratic extension $\dd$. In this section we will prove, that the existence of a $\G$-structure on $(\fg,\theta,\ip)$ implies $\dim \fl_-=3$ provided that the symmetric triple $(\fg,\theta,\ip)$ is indecomposable.
\begin{lm}\label{P3}
If $(\fg,\theta,\ip)$ is indecomposable and has a $\G$-structure, then $\dim \fl_-<3$ implies $\dim\fg_+\le3$.
\end{lm}
\proof If $(\fg,\theta,\ip)$ admits a $\G$-structure, then $\fg$ is solvable. By Lie's theorem $[\fg,\fg]$ acts on $\fg$ by nilpotent endomorphisms. Thus $\fg_+\subset [\fg,\fg]$ acts nilpotently on $\fg_-$. Hence there exists an element $U\not=0$ in $\fg_-$ such that $[\fg_+,U]=0$. Since $\fg$ is indecomposable $U$ is isotropic. 
Now let $({\cal O},[\psi])$ be a $\G$-structure. Then $U\cdot\psi\in\Delta_{\fg_-}$ is a $\fg_+$-invariant isotropic spinor and $\psi\perp U\cdot \psi$.
Now Corollary \ref{Ciso2} implies that the maximal dimension of an abelian Lie algebra contained in $\fg_+$ is three.

Finally, we will prove that $\fg_+$ is abelian, which will imply the assertion of the lemma. The metric Lie algebra $\fg_+$ is a quadratic extension of $\fl_+$ by $\fa_+$. In our situation $\dim \fl_+\le1$. Moreover, $\rho|_{\fl_+}=0$ by Lie's theorem since $\fl$ is solvable, $\fl_+=[\fl_-,\fl_-]\subset [\fl,\fl]$ and $\rho$ is semisimple. Thus $\fg_+$ is abelian.
\qed

\begin{pr}
If $(\fg,\theta,\ip)$ is indecomposable and admits a $\G$-structure, then $\dim \fl_-=3$.
\end{pr}
\proof If $\dim \fl_-<3$, then $\dim \fl_+\le1$. Let us first consider the case 
$\dim \fl_+=0$.
In this case $\fl$ is abelian, hence $\fl=\RR$ or $\fl=\RR^2$. If $\fl=\RR$, then $\alpha=\gamma=0$. Thus $\fa=\rho(\fl)(\fa)$ since $\fg$ is indecomposable. Because of $\dim \fa_-=7-2\cdot\dim \fl_-=5$ this shows that $\fa_+$ has also dimension 5. Hence $\fg_+=\fa_+$ is a five-dimensional abelian Lie algebra, which is a contradiction to Lemma~\ref{P3}. 

If $\fl=\RR^2=\Span\{Y,Z\}$, then $\fa_-$ has signature $(2,1)$. First consider the case $\alpha=0$. Every non-trivial indecomposable orthogonal $\RR^2$-module $(\bar\fa,\bar\rho)$ is of one of the following types
\benur
\item $\bar\fa\in\{\RR^{2,0},\RR^{2}\}$, $\bar \rho_{\Bbb C}$ has weights $\pm i\lambda$, $\lambda\in \fl^*$,
\item $\bar\fa=\RR^{1,1}$, $\bar \rho_{\Bbb C}$ has weights $\pm \lambda$, $\lambda\in \fl^*$,
\item $\bar\fa=\RR^{2,2}$, where $\fa_+=\RR^{1,1}$, $\fa_-=\RR^{1,1}$ and $\bar \rho_{\Bbb C}$ has weights $\pm \mu\pm i\nu$, $\mu,\nu\in \fl^*$.
\end{enumerate}

Since $\rho$ does not contain a trivial subrepresentation and $t_\rho=0$ by Prop.~\ref{Pricl} we have the following possibilities for  $(\fa,\rho)$:

(a)\  $\fa$ is the sum of two representations of type (i) with weights $\pm i\lambda_1$ and $\pm i \lambda_2$ and one representation of type (ii) with weights $\pm\lambda_3$. Then $t_\rho=0$ gives  
\begin{equation}\label{Elambda}
 -\lambda_1(L)\lambda_1(L')-\lambda_2(L)\lambda_2(L')+ \lambda_3(L)\lambda_3(L')=0\,.
 \end{equation}
Put $\lambda:=(\lambda_1,\lambda_2,\lambda_3)$. Then (\ref{Elambda}) says that $\lambda(Y)$ and $\lambda(Z)$ span an isotropic subspace of a pseudo-Euclidean space of signature $(2,1)$. The maximal dimension of such a subspace is one. Thus $\lambda(Y)$ and $\lambda(Z)$ are linearly dependent, which contradicts the indecomposability of $\fg$.

(b)\ $\fa$ is the sum of two representations of type (ii) and one representation of type (i). Then the argumentation is as in case (a). 

(c)\ $\fa$ is the sum of a representation of type (iii) and a representation of type (i) or (ii). Then 
$$\pm \lambda(L)\lambda(L')+4\mu(L)\mu(L')-4\nu(L)\nu(L')=0$$
holds for all $L,L'\in\fl$. As in (a) we see that $(\lambda(Y),\mu(Y),\nu(Y))$ and $(\lambda(Z),\mu(Z),\nu(Z))$ are linearly dependent, which contradicts indecomposability.

Now we consider the case $\alpha\not=0$. Then the $(\fl,\theta_\fl)$-module $\fa$ decomposes into the one-dimensional image of $\alpha$ and a six-dimensional $(\fl,\theta_\fl)$-module $\fa'$, which is equivalent to the module $\fa$ considered in the case $\alpha=0$. Hence $\fg_+=\fa_+$ is four-dimensional. Thus $\fg$ cannot obtain a $\G$-structure by Lemma~\ref{P3}. 

Now assume $\dim\fl_+=1$.
Then $\dim\fl_-=2$ and $\fl\in\{\fh(1),\, \fn(2),\, \fr_{3,-1} \}$, see \cite{KO2}, Prop.~7.2.
Moreover, $\fa_-$ has signature $(2,1)$. 
Since $\rho$ is semisimple and $\fl$ is solvable $\rho$ is a representation of the abelian Lie algebra $\fl/[\fl,\fl]$. We decompose $\fa=\fa^\fl\oplus\fa'$ with $\fa'=\rho(\fl)(\fa)$. Since $\fa^\fl\subset \alpha(\fl,\fl)$ we have $\dim \fa^\fl_-\le 2$. Thus $\dim \fa'_-\ge1$. If $\dim\fa'_+=\dim\fa'_->1$, then we are done. Hence assume $\dim \fa^\fl_-= 2$ and $\dim \fa'_+=\dim\fa'_-=1$. The latter equation implies that  $\fl\not=\fh(1)$ since the Killing form of $\fh(1)$ vanishes, which would imply  $t_\rho=0$ by Prop.~\ref{Pricl}, thus $\rho=0$, a contradiction.
Furthermore, $\dim \fa^\fl_-= 2$ implies that $\fl\not\in=\{\fn(2),\ \fr_{3,-1}\}$, since in both cases $\fa^\fl_-$ is spanned by $\alpha(Y,Z)$ for any $[\alpha,\gamma]\in \HQtheta_0$, see \cite{KO2}, Prop.~7.3 (note that there is an obvious typo). 
\qed
\subsection{The canonical isotropic ideal and the cross product} 
\label{S42}
Let $(\fg,\theta,\omega,\ip)$ be a symmetric triple with $\G$-structure and let $b$ be the cross product that corresponds to $\omega$. We will prove that the canonical isotropic ideal $\fri_-=\fri\cap \fg_-$ is invariant under $b$. To this end we consider the $\G$-structure $\omega$ equivalently as a pair $({\cal O},[\psi])$ according to the remark after Definition~\ref{Dg2}. Furthermore, let $\ph\in\Delta_{\fg_-}$, $\ph\not=0$, be an isotropic spinor corresponding to $U_\ph=\fri_-$, i.e., $\fri_-\cdot\ph=0$. Recall that $\ph$ is uniquely determined up to multiplication with a real number $r\not=0$.
\begin{lm}
The spinor $\ph$ is $\fg_+$-invariant. 
\end{lm}
\proof Let $X\in\fg_+$ be arbitrary. Note that $Y\cdot A \cdot \ph= A\cdot Y\cdot\ph- \lambda_*(A)(Y)\cdot\ph$ holds for all $Y\in\fg_-$ and $A\in\fspin(\fg_-)$. In particular, take $Y\in U_\ph$. Then
$$Y\cdot \tad(X)\cdot\ph=\tad(X)\cdot Y\cdot\ph -[X,Y]\cdot\ph=-[X,Y]\cdot\ph.$$ 
Since $U_\ph=\fri_-$ is $\fg_+$-invariant, we have $[X,Y]\in U_\ph$ and thus $Y\cdot \tad(X)\cdot\ph=0$ for all $Y\in U_\ph$. Hence $\tad(X)\cdot\ph=t\ph$ for some $t\in\RR$. 
On the other hand, $X\in\fg_+\subset[\fg,\fg]$ acts nilpotently on $\Delta_{\fg_-}$ since $\fg$ is solvable. Hence $t=0$, which proves the claim.
\qed
\begin{lm}
If $\fg$ is indecomposable, then $\ph\perp\psi$.
\end{lm}
\proof Assume that $\la\ph,\psi\ra_\Delta\not=0$. Then we can choose $\fg_+$-invariant orthogonal elements $\psi_1,\psi_2\in\Delta_{\fg_-}$ such that  $\la\psi_i,\psi_i\ra_\Delta\not=0$, $i=1,2$. By Proposition \ref{Pmult} we can define a vector $X$ in $\fg_-$ such that $\psi_2=X\cdot \psi_1$. Then $X$ is not isotropic and satisfies $[\fg_+,X]=0$. This is a contradiction to the indecomposability of $\fg$. \qed
\begin{pr}\label{Pfri} The cross product $b$ has the following properties
\begin{enumerate}[label=(\roman*)]
\item $b(\fri_-^\perp,\fri_-^\perp)\subset \fri_-$,   
\item $\fn^*:=b(\fri_-,\fri_-)$ satisfies  $\fn^*\not=0$, $b(\fn^*,\fri_-)=0$ and 
$b(\fn^*,\fri_-^\perp)=\fn^*$. 
\end{enumerate}
\end{pr}
\proof Any oriented (with respect to the orientation induced by $\omega$) Witt basis ${\bf b}=(b_1,\dots,b_7)$ of $\fg_-$ gives us an isometry from $\fg_-$ to $\RR^{4,3}$. Thus we can identify also $\Delta_{\fg_-}\cong\Delta_{4,3}$. By Prop.~\ref{Ciso2} we can choose ${\bf b}$ in such a way that $\psi=s_1+s_5$ and $\ph=s_6$. Then $\fri_-=U_\ph=\Span\{b_1,b_2,b_3\}$ and $\fri_-^\perp=\Span\{b_1,\dots,b_4\}$. Furthermore, $\omega$ is equal to the 3-form $\omega_0$, which was defined by Equation~(\ref{Eomega}). Hence 
\begin{eqnarray*}
&b(b_1,b_2)=\sqrt2 b_3,\ b(b_1,b_3)=b(b_2,b_3)=0&\\
&b(b_1,b_4)=-b_1,\ b(b_2,b_4)=-b_2,\ b(b_3,b_4)=b_3&
\end{eqnarray*}
thus $\fn^*=\RR\cdot b_3$ and the assertion follows.
\qed

Using Prop.~\ref{Pfri}, (ii) we can define an orientation of the one-dimensional vector space $\fri_-^\perp/\fri_-$ in the following way. Note first that $b:\fri_-^\perp/\fri_-\otimes\fn^*\rightarrow \fn^*$ correctly defined. 
\begin{de}
 A vector $A$ in $\fri_-^\perp/\fri_-$ is said to have a positive orientation if $b(A,U)$ is a positive multiple of $U$ for all $U\in\fn^*$. 
\end{de}
\begin{de} A Lie algebra with $B$-structure $(\fl,\theta_\fl,b_\fm)$ consists of a solvable Lie algebra with involution $(\fl,\theta_\fl)$ satisfying $[\fl_-,\fl_-]=\fl_+$  and a non-trivial antisymmetric bilinear map $b_\fm:\fm\times\fm\rightarrow\fl_-$ defined on a two-dimensional subspace $\fm$ of $\fl_-$ such that 
\begin{itemize}
\item[{\rm (L1)}] $[\fl_+,\fm]\subset \fm$,
\item[{\rm (L2)}] $[\fm,\fm]_\fl=0$,
\item[{\rm (L3)}] $\fn:=b_\fm(\fm,\fm)$ is complementary to $\fm$ in $\fl_-$.
\end{itemize} 
An isomorphism of Lie algebras with $B$-structure $S:(\fl,\theta_\fl,b_\fm)\rightarrow(\fl',\theta_{\fl'},b_{\fm'})$ is an isomorphism of Lie algebras with involutions that satisfies in addition $S(\fm)=\fm'$ and 
$$b_{\fm'}(SL,SL')=Sb_\fm(L,L')$$ holds for all $L,L'\in\fm$.
\end{de}
\begin{re}\label{Riii}{ \rm  Since $\fl$ is solvable and $\fl_+\subset[\fl,\fl]$ the operator $\ad(L_+)$ is nilpotent for all $L_+\in\fl_+$. 
Hence (L1) implies 
\begin{equation}\label{Elm}
[\fl_+,\fl_-]\subset\fm
\end{equation}
and  $\tr \ad (L_+)|_\fm=0$ for all $L_+\in\fl_+$. The latter equation is equivalent to
\begin{equation}\label{EOh2}
b_\fm([L_+,L_1]_\fl,L_2)+b_\fm(L_1,[L_+,L_2]_\fl)=0.
\end{equation}
Indeed, the left hand side of (\ref{EOh2}) equals $\tr \ad(L_+)|_\fm \cdot b_\fm(L_1,L_2)$.
}\end{re}
Let $(\fa,\ip_\fa,A)$ be a pseudo-Euclidean space with distinguished time-like unit vector $A$, i.e., $\la A,A\ra_\fa=-1$. We define an involution $\theta_\fa$ on $\fa$ by 
\begin{equation}\label{Einv}
\fa_-=\RR\cdot A,\quad \fa_+=\fa_-^\perp.
\end{equation}
We will call $\theta_\fa$ the induced involution. The vector $A$ may be considered as an orientation of $\fa_-$.

\section{The standard model}
\subsection{Quadratic cocycles of Lie algebras with $B$-structure}
Let $(\fl,\theta_\fl,b_\fm)$ be a Lie algebra with $B$-structure and $(\fa,\ip_\fa\,A)$ a pseudo-Euclidean space with distinguished time-like unit vector $A$. Let $\theta_\fa$ be the induced involution as defined in (\ref{Einv}). We consider $\fa$ as a trivial $(\fl,\theta_\fl)$-module.
\begin{de}
 We define $\ZQb$ to be the set of those $(\alpha,\gamma)\in\ZQtheta$ that satisfy
\begin{itemize}
\item[{\rm (Z1)}] $\alpha(\fm,\fm)=0$;
\item[{\rm (Z2)}]$[L,b_\fm(L',L'')]=\la \alpha(L,L''), A\ra L'-\langle \alpha(L,L'),A\ra L'' $ for all $L\in\fl_+, L',L''\in \fm$;
\item[{\rm (Z3)}] $2 \gamma(L,L',L'') =-\langle A,\alpha(L,b_\fm(L',L''))\rangle$ for all $L\in\fl_+$, $L',L''\in \fm$.
\end{itemize}
We will say that an element of $\ZQb$ is admissible if it is admissible as an element  of $\ZQtheta$ and we will denote by $\ZQb_\sharp$ the set of admissible elements of $\ZQb$.
\end{de}
In the following we will need the group
\begin{equation}\label{EcN}
{\cal N}:=\{S_0\in\Aut(\fl,\theta_\fl)\mid S_0(\fm)\subset\fm,\ S_0|_\fm=\Id_\fm,\ \bar S_0=\Id_{\fl_-/\fm},\ S_0|_{\fl_+}=\Id\},
\end{equation}
where $\bar S_0$ denotes the map induced by $S_0$ on $\fl_-/\fm$.
\begin{de}
We will say that elements $(\alpha,\gamma)$ and $(\alpha',\gamma')$ of $\ZQb$ are equivalent if and only if there exist an isomorphism $S_0\in{\cal N}$ and an element $(\tau,\sigma)\in\CQtheta$ such that 
\benur
\item $(\tau,\sigma)$ has the properties 
\begin{itemize}
\item[{\rm (B1)}]$\tau(L)=\tr (\proj_\fm S_0(b_\fm(L,\cdot)))\cdot A$ for all $L\in\fm$,\\
where $\proj_\fm:\fl=\fm\oplus\fn\rightarrow\fm$ is the projection,
\item[{\rm (B2)}] $ 2\sigma(L',L'')=-\la\tau(b_\fm(L',L'')),A\ra$ for all $L',L''\in\fm;$
\end{itemize}
\item $(\alpha,\gamma)=(S_0^*\alpha', S_0^*\gamma')\cdot (\tau,\sigma).$
\end{enumerate}
We will denote the set of equivalence classes by $\HQb$.
\end{de}
\subsection{Construction of the standard model}
Let $(\fl,\theta_\fl,b_\fm)$ be a Lie algebra with $B$-structure. Recall that, in particular, $(\fl,\theta_\fl)$ is proper. Let $(\fa,\ip_\fa\,A)$ be a pseudo-Euclidean space with distinguished time-like unit vector $A$. Let $\theta_\fa$ be the induced involution as defined in (\ref{Einv}). As above we consider $\fa$ as a trivial $(\fl,\theta_\fl)$-module. Let $(\alpha,\gamma)$ belong to $\ZQb_\sharp$. Now we consider the quadratic extension $\fd:=\dd$ of $(\fl,\theta_\fl)$ by $(\fa,\ip_\fa,\theta_\fa)$ corresponding to $(\alpha,\gamma)\in\ZQtheta$. 
Let $L_1,L_2$ be a basis of $\fm$ and $L_3:=(1/\sqrt2)\cdot b_\fm(L_1,L_2)$. Furthermore, let $Z_1,Z_2,Z_3\in \fl^*$ be the dual basis of $L_1,L_2,L_3$. Let $\omega$ be the 3-form on $\fg_-$ that is defined by Equation~(\ref{Enochmal}) with respect to the dual basis $\sigma^1,\dots,\sigma^7$ of $Z_1,Z_2,Z_3,A,L_1,L_2,L_3$.
It can be checked easily that $\omega$ does not depend on the choice of the basis $L_1,L_2$ of $\fm$. By definition $\omega$ is a nice 3-form. 
\begin{pr}\label{PS} Under the above assumptions $\omega$ is invariant under $\fg_+$.
\end{pr}
\proof  The computations become simpler when we use the spinorial description of the $\G$-structure. The orientation defined by $\omega$ is such that $Z_1,Z_2,Z_3,A,L_1,L_2,L_3$ is positively oriented. 
We can define a representation $\phi :\Cl(\fg_-)\rightarrow \End(\Delta_{\fg_-})$ of type one by $\phi(Z_i)=\Phi(e_i)$, $\phi(A)=\Phi(e_4)$ and $\phi(L_i)=\Phi(e_{i+4})$, $i=1,2,3$. Then $\omega$ corresponds to  $\psi= s_1+s_5$ according to Proposition \ref{Pobp}. Hence it suffices to show that $\psi$ is $\fg_+$-invariant. From (\ref{Etilde}) and (\ref{EPhi}) we obtain
\begin{eqnarray}
4\,\tad(Z_+)\cdot\psi&=&\sum_{i=1}^3 Z_i\cdot [Z_+,L_i]\cdot\psi =2\sum_{1\le i<j\le3}Z_+([L_i,L_j])Z_iZ_j\cdot\psi \nonumber\\
&=&2\sum_{1\le i<j\le3}Z_+([L_i,L_j])e_i e_j\cdot(s_1+s_5)=4Z_+([L_1,L_2])s_6, \label{EadZ}
\end{eqnarray}
which vanishes by assumption. In the same way we get
\begin{eqnarray*}
4\,\tad(A_+)\cdot\psi&=&\sum_{i=1}^3 Z_i\cdot [A_+,L_i]\cdot\psi =2\sum_{1\le i<j\le3}\la A_+,\alpha(L_i,L_j)\ra Z_iZ_j\cdot\psi\\
&=& 2\sum_{1\le i<j\le3}\la A_+,\alpha(L_i,L_j)\ra e_ie_j\cdot(s_1+s_5)=4\la A_+,\alpha(L_1,L_2)\ra s_6=0.
\end{eqnarray*}
Furthermore,
\begin{eqnarray}\label{EadL}
4\,\tad(L_+)\cdot\psi&=&\Big(\sum_{i=1}^3\Big(Z_i\cdot[L_+,L_i]+L_i\cdot[L_+,Z_i]\Big)-A\cdot[L_+,A]\Big)\cdot\psi \nonumber\\
&=& \sum_{i=1}^3\Big(Z_i\cdot([L_+,L_i]_\fl+\alpha(L_+,L_i)+\gamma(L_+,L_i,\cdot))-L_i\cdot Z_i([L_+,\cdot]_\fl)\Big)\cdot\psi \nonumber\\
&& +A\cdot\la A,\alpha(L_+,\cdot)\ra \cdot\psi \nonumber\\
&=& 2\sum_{i,j=1}^3Z_j([L_+,L_i])Z_iL_j\cdot\psi+2\sum_{1\le i<j\le3} \gamma(L_+,L_i,L_j)Z_iZ_j\cdot\psi \nonumber\\
&&-2\sum_{k=1}^3 \la A,\alpha(L_+,L_i)\ra Z_iA\cdot\psi,
\end{eqnarray}
where we used $\tr (\ad L_+)|_{\fl_-}=0$. We will show that the right hand side of this equation vanishes.
We have
\begin{eqnarray}\label{EadLs}
\lefteqn{\Big(\sum_{i,j=1}^3Z_j([L_+,L_i])e_ie_{j+4}+\sum_{i<j} \gamma(L_+,L_i,L_j)e_ie_j-\sum_{k=1}^3 \la A,\alpha(L_+,L_i)\ra e_ie_4\Big)\cdot(s_1+s_5)}\nonumber \\
&=&2\Big( -Z_1([L_+,L_1])s_5-Z_1([L_+,L_3])s_3-Z_2([L_+,L_2])s_5-Z_2([L_+,L_3])s_8+\qquad \nonumber\\
&&Z_3([L_+,L_1])s_7+Z_3([L_+,L_2])s_4-Z_3([L_+,L_3])s_1\Big)+2\gamma(L_+,L_1,L_2) s_6 +\nonumber\\
&&\sqrt2\Big(- \la A,\alpha(L_+, L_1)\ra s_8 + \la A,\alpha(L_+, L_2)\ra s_3+ \la A,\alpha(L_+, L_3)\ra s_6\Big).
\end{eqnarray}
Now we use that
$$  -Z_1([L_+,L_1])-Z_2([L_+,L_2])=-\tr \ad(L_+)|_\fm=0,$$
see Remark~\ref{Riii}. Furthermore,
$$2Z_1([L_+,L_3])=\sqrt2 \la A,\alpha(L_+, L_2)\ra,\quad 2Z_2([L_+,L_3])= -\sqrt2 \la A,\alpha(L_+, L_1)\ra$$
since $\sqrt 2[L_+,L_3]=[L_+,b_\fm(L_1,L_2)]=\la \alpha(L_+,L_2), A\ra L_1-\langle \alpha(L_+,L_1),A\ra L_2 $ by (${\rm Z2}$), and
$$Z_3([L_+,L_1])=Z_3([L_+,L_2])=Z_3([L_+,L_3])=0$$
since $[\fl_+,\fl_-]\subset\fm$. Finally,
$$2\gamma(L_+,L_1,L_2)+ \sqrt2\la A,\alpha(L_+, L_3)\ra=2\gamma(L_+,L_1,L_2)+ \la A,\alpha(L_+, b_\fm(L_1,L_2))\ra=0$$
by {\rm (Z3)}.
\qed
\begin{co}
If $(\fl,\theta_\fl,b_\fm)$ is a Lie algebra with $B$-structure, $(\fa,\ip_\fa,A)$ a pseudo-Euclidean space with distinguished time-like unit vector and if $(\alpha,\gamma)$ is in $\ZQb_\sharp$, then $\ddm:=(\fd, \theta,\omega,\ip)$ is a symmetric triple with $\G$-structure.  
\end{co}
We consider the Lie algebras
\begin{eqnarray*}
\fh(1)&=&\{[X,Y]=Z\}\\
\fg_{4,1}&=&\{B=[L_2,L_3],\ [B,L_3]=L_1\},
\end{eqnarray*}
where we use the following convention. The vectors appearing on the right hand side constitute a basis of the Lie algebra and all brackets of basis vectors not mentioned are equal to zero. The Lie algebra $\fh(1)$ is the three-dimensional Heisenberg algebra, $\fg_{4,1}$ is the only indecomposable real nilpotent Lie algebra of dimension 4.
\begin{pr}\label{Pl}
Each indecomposable symmetric triple with $\G$-structure  is isomorphic to a quadratic extension $\ddm$ for some $(\alpha,\gamma)\in\ZQb_\sharp$, where $(\fl,\theta_\fl,b_\fm)$ is isomorphic to one of the following two triples:

$\fl_-=\Span\{L_1,L_2,L_3\}$, $\fm=\Span\{L_1,L_2\}$, $b_\fm(L_1,L_2)=\sqrt 2 L_3$, \\
$\fl_+=\Span\{B=[L_2,L_3]_\fl\}$, 
 $B\not=0$,  and
\begin{enumerate}
\item $[B,L_1]=[B,L_2]_\fl=0,\ [B,L_3]_\fl=L_1$, or
\item $\ad(B)=0$.
\end{enumerate}
In the first case $\fl$ is isomorphic to $\fg_{4,1}$, in the second one to $\fh(1)\oplus\RR$.
\end{pr}
\proof Let $(\fg,\theta,\omega,\ip)$ be an indecomposable symmetric triple with $\G$-structure. Let $\omega$ correspond to the pair $({\cal O},[\psi])$ according to the remark after Definition~\ref{Dg2}. In Section \ref{S32} we have seen that the intersection $\fri_-$ of the canonical isotropic ideal $\fri$ of $\fg$ with $\fg_-$ is three-dimensional. We proceed as in Section~\ref{S42}, especially as in the proof of Proposition~\ref{Pfri}. Let $\ph\not=0$ be an isotropic spinor that satisfies $\fri_-\cdot\ph=0$. Any choice of a Witt basis ${\bf b}$ of $\fg_-$ gives an equivalence $\Delta_{\fg_-}\cong \Delta_{4,3}$. By Proposition~\ref{Ciso2} we can choose ${\bf b}$ such that $\psi=s_1+s_5$ and $\ph=s_6$.
In particular, 
$\fri_-=\Span\{b_1,b_2,b_3\},\ \fri_-^\perp=\{b_1,b_2,b_3,b_4\}$, 
and $\omega$ is given by (\ref{Eomega}). We choose an  isotropic complement
$V_+\subset\fg_+$ of $\fri_+^\perp$ in $\fg_+$, and we put
$V_-:=\Span\{b_5,b_6,b_7\}$ and $V:=V_+\oplus V_-$.
Then we may identify $V$ and $\fl:=\fg/\fri^\perp$ as vector spaces with involution. Moreover, we put $\fa:=(\fri\oplus V)^\perp$ and $A:=b_4$. Then $\fa$ is a pseudo-Euclidean space with distinguished time-like unit vector $A$ and the induced involution $\theta_\fa$ coincides with $\theta|\fa$. In particular, the symmetric triple $(\fg,\theta,\ip)$ is an admissible quadratic extension $\dd$ of $(\fl,\theta_\fl)$ by $\fa$.  

Note that $b(V_-, V_-)\subset V_-$. More exactly,
$$b(b_5,b_6)=\sqrt2 b_7 ,\quad b(b_5,b_7)=b(b_6,b_7)=0. $$
We put 
$\fm:=\Span\{b_5,b_6\}$ and $b_\fm= b|_{\fm\times \fm}$. We want to show that $(\fl,\theta_\fl,b_\fm)$ is a Lie algebra with $B$-structure, that the representation $\rho$ of $\fl$ on $\fa$ is trivial and that $(\alpha,\gamma)$ is in $\ZQb$. 

The existence of $\omega$ implies that $\fg$ is solvable, see Proposition~\ref{Pric}. Hence, also $\fl$ is solvable. Since $(\fg,\theta,\ip)$ is a symmetric triple $[\fg_-,\fg_-]=\fg_+$ holds, which implies $[\fl_-,\fl_-]=\fl_+$. We have to show that $\rho$ is trivial and to verify (L2), (L3) and (Z1) -- (Z3). As for $\rho$, we already know that $\rho|_{\fl_+}=0$ since $\fl$ is solvable and $\rho$ is semisimple.
Even if $\rho\not=0$ Equation (\ref{EadZ}) remains true, which implies (L2).

Furthermore,
\begin{eqnarray*}
4\,\tad(A_+)\cdot\psi&=&\sum_{i=1}^3 Z_i\cdot [A_+,L_i]\cdot\psi -A\cdot[A_+,A]\cdot \psi \\
&=&\Big( 2\sum_{i=1}^3 \la \rho(L_i)(A_+),A\ra Z_i  A +2\sum_{i<j}\la A_+,\alpha(L_i,L_j)\ra Z_iZ_j\Big) \cdot\psi\\
&=& \Big(2\sum_{i=1}^3 \la \rho(L_i)(A_+),A\ra e_i  e_4+2\sum_{i<j}\la A_+,\alpha(L_i,L_j)\ra e_ie_j\Big)\cdot(s_1+s_5)\\
&=& 2\sqrt2(\la \rho(L_1)(A_+),A\ra s_8-\la \rho(L_2)(A_+),A\ra s_3-\la \rho(L_3)(A_+),A\ra s_6)\\
&&+4\la A_+,\alpha(L_1,L_2)\ra s_6=0
\end{eqnarray*}
implies $\rho|_\fm=0$ and
\begin{equation}\label{E16}
 \rho(L_3)(A)+\sqrt 2\alpha(L_1,L_2)=0.
\end{equation}
Since $\rho|_{\fl_+}=0$ Equation (\ref{EadL}) remains true. Hence (\ref{EadLs}) holds and implies $(L1)$, $(Z2)$ and $(Z3)$. It remains to prove $\alpha(\fm,\fm)=0$ and $\rho|_\fn=0$ and to determine all possible $(\fl,\theta_\fl,b_\fm)$.

Let $(\fl,\theta_\fl,b_\fm)$ be a Lie algebra with $B$-structure and let $(\rho,\fa)$ be an orthogonal $(\fl,\theta_\fl)$-module with $\rho|_\fm=0$. Assume that $(\alpha,\gamma)\in\ZQtheta$ is admissible and satisfies (Z2), (Z3) and (\ref{E16}). Let $L_1,L_2$ be a basis of $\fm$. Then
$B:=[L_2,L_3]_\fl$ and $C:=[L_3,L_1]_\fl$ span $\fl_+$.
The Jacobi identity for $L_1,L_2,L_3$ implies 
\begin{equation}\label{Ebald}
[B,L_1]+[C,L_2]=0.
\end{equation} 
Since $B$ and $C$ act nilpotently on $\fl_-$, we have $\tr \ad(B)|_\fm=\tr \ad(C)|_\fm=0$ by (\ref{Elm}).  Together with (\ref{Ebald}) this gives
$$\begin{array}{ll}
[B,L_1]_\fl =aL_1+dL_2,&  [B,L_2]_\fl=bL_1-aL_2,\\{} [C,L_1]_\fl=dL_1+cL_2,& [C,L_2]_\fl=-aL_1-dL_2.
\end{array}$$
Because of $[\fl_+,\fl_-]\subset\fm$ and $[\fm,\fm]=0$ the Jacobi identity for $B,L_2,L_3$ and for $C,L_1,L_3$ gives $[[B,L_2]_\fl,L_3]_\fl=0$ and $[[C,L_1]_\fl,L_3]_\fl=0$. 

First we consider the case that $B$ and $C$ are linearly independent. Then the equation $[[B,L_2]_\fl,L_3]_\fl=0$ yields $a=b=0$ and $[[C,L_1]_\fl,L_3]_\fl=0$ gives $d=c=0$. Hence $\ad(L)|_{\fm}=0$ for all $L\in\fl_+$.
Furthermore, the Jacobi identity for $B,L_1,L_3$ gives $[B,C]=0$.
Put
$$r:=\la \alpha(B,L_2),A\ra,\ s:= -\la\alpha(B,L_1),A\ra, \ p:=\la\alpha(C,L_2),A\ra,\ q:=-\la\alpha(C,L_1),A\ra$$
Condition (${\rm Z2}$)  yields
\begin{equation}\label{EBL}
\sqrt 2[B,L_3]_\fl=rL_1+sL_2,\quad \sqrt 2[C,L_3]_\fl=pL_1+qL_2.
\end{equation}
Assume that $\rho(L_3)(A)=0$. Then $\alpha(\fm,\fm)=0$ by (\ref{E16}), thus
$$0=d\alpha(B,L_1,L_3)=-\alpha([L_1,L_3],B)-\alpha([L_3,B],L_1)-\alpha([B,L_1],L_3)=-\alpha(B,C).$$
In particular, $B,C,\sigma_B,\sigma_C$ span an abelian subalgebra of $\fg_+$, which is impossible by Proposition~\ref{Ciso2}. Hence $\rho(L_3)(A)=:A_0\not=0$. Then $\sqrt2 \alpha(L_1,L_2)=-A_0$. Now
\begin{eqnarray*}
d\alpha(B,L_2,L_3)&=&\rho(L_3)(\alpha(B,L_2))-\alpha([L_3,B],L_2)\nonumber \\
&=& - r \rho(L_3)(A) +(1/\sqrt2) \alpha(rL_1+sL_2,L_2)\nonumber \\
&=&  (-r-r/2)A_0\, =0
\end{eqnarray*}
and, analogously,
$$d\alpha(C,L_1,L_3)= (q+q/2)A_0 =0,$$
which gives $q=r=0$. Using this we get
\begin{eqnarray}
d\alpha(L_1,L_2,L_3)&=&\rho(L_3)(\alpha(L_1,L_2))-\alpha([L_2,L_3],L_1)-\alpha([L_3,L_1],L_2)\nonumber \\
&=& - (1/\sqrt2) \rho(L_3)(A_0) -\alpha(B,L_1)-\alpha(C,L_2)\nonumber \\
&=& - (1/\sqrt2)\rho(L_3)(A_0) -sA+pA\ =0\label{Ea}\\[2ex]
d\alpha(B,C,L_3) &=& \rho(L_3)( \alpha(B,C))-\alpha([C,L_3],B)-\alpha([L_3,B],C)\nonumber \\
&=& \rho(L_3)(\alpha(B,C))-(p/\sqrt2)\alpha(L_1,B)+(s/\sqrt2)\alpha(L_2,C)\nonumber\\
&=& \rho(L_3)(\alpha(B,C))+\sqrt2psA\label{Eb}
\\[2ex]
d\alpha(B,L_1,L_3) &=& \rho(L_3)(\alpha(B,L_1))-\alpha([L_1,L_3],B)-\alpha([L_3,B],L_1)\nonumber \\
&=&s \rho(L_3)(A)-\alpha(B,C)+(s/\sqrt2)\alpha(L_2,L_1)\nonumber\\
&=& -\alpha(B,C)+(s+s/2)A_0\label{Ec}
\end{eqnarray}
and, analogously,
\begin{equation}
d\alpha(C,L_2,L_3) =-\alpha(B,C)+(-p-p/2)A_0. \label{Ed}
\end{equation}
Equations (\ref{Ec}) and (\ref{Ed}) give $s=-p$. Note that $s\not=0$ since as above $\alpha(B,C)\not=0$ holds by Prop.~\ref{Ciso2}. By (\ref{Eb}) and (\ref{Ec}) we get $\rho(L_3)(A_0)=(2\sqrt 2 s/3)A$.
On the other hand, $\rho(L_3)(A_0)=-2\sqrt 2 sA$ by (\ref{Ea}), which gives a contradiction.

Now suppose $\dim \fl_+=1$. We may assume $[L_3,L_1]_\fl=0$ and $\fl_+=\RR\cdot B$ with  $B:=[L_2,L_3]_\fl$.
The Jacobi identity for $L_1,L_2,L_3$ gives $[B,L_1]_\fl=0$. Recall from Remark~\ref{Riii} that $\tr \ad(B)|_\fm=0$, thus $[B,L_2]_\fl=bL_1$ for some $b\in\RR$. Put 
$$r:=\la A,\alpha(B,L_2)\ra,\ s:=-\la A,\alpha(B,L_1)\ra.$$
Then Condition (${\rm Z2}$) gives $$\sqrt 2 [B,L_3]=rL_1+sL_2.$$
As above, $A_0:=\rho(L_3)(A)$. Computations analogous to (\ref{Ea}) and (\ref{Ec})  give
\begin{eqnarray*}
d\alpha(L_1,L_2,L_3)&=& - (1/\sqrt2) \rho(L_3)(A_0) - sA\ =0 \\
d\alpha(B,L_1,L_3)&=&  (s+s/2)A_0\, =0,
\end{eqnarray*}
which implies $A_0=0$ and $s=0$. In particular, $\alpha(B,L_1)=0$ and $\sqrt 2 [B,L_3]=rL_1$.

From $d\alpha=0$ we get $\alpha([B,L_2],L_3)=\alpha([B,L_3],L_2)=(r/\sqrt2)\alpha(L_1,L_2)=0$, thus $b\alpha(L_1,L_3)=0$. For $b\not=0$ we would get $\alpha(L_1,\cdot)=0$. On the other hand, in this case $\fl\cong\fg_{4,1}$ and for $\fg_{4,1}$ the admissibility condition ($A_2$) in \cite{KO2}, Definition~5.3. implies that a cocycle $(\alpha,\gamma)$ can only be admissible if $\alpha(\fz,\cdot)=\alpha(L_1,\cdot)\not=0$,  see also \cite{Knil}, Prop.~4.4 for this fact. Hence $b=0$.
Rescaling $L_1,L_2$ and $L_3$ we may assume $r\in\{0,1\}$.

If $B=C=0$, then the Killing form of $\fl$ vanishes, thus $\rho=0$ by Proposition~\ref{Pricl}. Since $(\fg,\theta,\ip)$ is indecomposable this yields $\fa_-=\alpha(\fl_+,\fl_-)=0$, which contradicts $\dim \fa_-=1$.
\qed

\section{Classification}
\subsection{Classification in terms of $\HQb$}
In this subsection we want to describe the isomorphism classes of quadratic extensions $\ddm$ in terms of $\HQb$.

For a Lie algebra with $B$-structure $(\fl,\theta_\fl,b_\fm)$ we define
$${\cal G}:=\left\{\ S\ \left| \begin{array}{l}S\in \Aut(\fl,\theta_\fl),\ S(\fm)=\fm,\\[1ex]\proj_\fn(S(b_\fm(L,L'))=b_\fm(S(L),S(L'))\mbox{ for all } L,L'\in\fm \end{array}\right.\right\}.$$
Then ${\cal G}$ equals the semi-direct product $\Aut(\fl,\theta_\fl,b_\fm)\ltimes {\cal N}$, where $\cN$ was defined in (\ref{EcN}) and 
$$\Aut(\fl,\theta_\fl,b_\fm)=\{S\in {\cal G}\mid S(\fn)=\fn\}$$
is the automorphism group of $(\fl,\theta_\fl,b_\fm)$.

\begin{de}\label{Dae} We will say that $\fd_{\alpha,\gamma}(\fl,\theta_\fl,b_\fm,\fa)$ and $\fd_{\alpha',\gamma'}(\fl,\theta_{\fl},b_{\fm},\fa)$ are equivalent if there is an isomorphism $\Psi_0: \fd_{\alpha,\gamma}(\fl,\theta_\fl,b_\fm,\fa) \rightarrow \fd_{\alpha',\gamma'}(\fl,\theta_{\fl},b_{\fm},\fa)$ of symmetric triples with $\G$-structure such that 
\benur
\item $\Psi_0(\fl^*)=\fl^*$ (hence $\Psi_0(\fl^*\oplus\fa)=\fl^*\oplus\fa)$, 
\item $\Psi_0|_\fa=\Id \ \mod \fl^*$,
\item the map $S_0:=\proj_\fl\circ\Psi_0|_\fl:\fl\rightarrow \fl$ belongs to ${\cal N}$.
\end{enumerate}
\end{de}
\begin{pr}\
The quadratic extensions $\fd_{\alpha,\gamma}(\fl,\theta_\fl,b_\fm,\fa)$ and $\fd_{\alpha',\gamma'}(\fl,\theta_{\fl},b_{\fm},\fa)$ are equivalent if and only if $[\alpha,\gamma]=[\alpha',\gamma']\in\HQb$.
\end{pr}
\proof In \cite{KO2} we proved that the isomorphisms $\Psi_0: \fd_{\alpha,\gamma}(\fl,\theta_\fl,\fa) \rightarrow \fd_{\alpha',\gamma'}(\fl,\theta_{\fl},\fa)$ of symmetric triples that satisfy Properties (i) -- (iii) of Definition \ref{Dae} are exactly the linear maps 
\begin{equation}\label{EM1}
    \Psi_0=\left(
\begin{array}{ccc}
    (S_0^*)^{-1} & -(S_0^*)^{-1}\tau^{*} & (S_0^*)^{-1}(\bar\sigma -\frac 12\tau^{*}\tau ) \\
    0 & \Id & \tau  \\
    0 & 0 & S_0
\end{array}\right) :\fl^{*}\oplus\fa\oplus\fl \longrightarrow 
\fl^{*}\oplus\fa\oplus\fl,
\end{equation}
where $(\tau,\sigma)$ is in $\CQtheta$ for $\sigma(L_1,L_2):=\bar\sigma(L_1)(L_2)$, $S_0:\fl\rightarrow \fl$ is in $\Aut(\fl,\theta_\fl,\fm)$ and $(\alpha,\gamma)=(S_0^*\alpha', S_0^*\gamma')\cdot (\tau,\sigma).$ It remains to decide which of these maps preserve the $\G$-structure. Choose a basis $L_1,L_2$ of $\fm$ and put $L_3:=\sqrt 2 b_\fm(L_1,L_2)$. Hence $L_1,L_2,L_3$ is a basis of $\fl$. Denote by $Z_1,Z_2,Z_3$ the dual basis of $\fl^*$. With respect to these bases we have
$$S_0(L_3)=L_3+s_1 L_1+s_2 L_2,\quad \tau(L_i)= t_iA,\ i=1,2,3.$$
Since 
\begin{eqnarray*}
b(\Psi_0(Z_3),\Psi_0(L_3))&=&b(Z_3,L_3+s_1 L_1+s_2 L_2+t_3A)\\
&=&-A+s_1\sqrt2 Z_2-s_2\sqrt2 Z_1-t_3Z_3
\end{eqnarray*}
and 
$$\Psi_0(b(Z_3,L_3))=\Psi_0(-A)=-A-t_1Z_1-t_2Z_2-(t_3+s_1t_1+s_2t_2)Z_3$$
the equation $b(\Psi_0(Z_3),\Psi_0(L_3))=\Psi_0(b(Z_3,L_3))$ is equivalent to 
$$ t_1=\sqrt 2 s_2,\quad t_2=-\sqrt2 s_1,$$
hence to (B1).
Furthermore,
\begin{eqnarray*}
b(\Psi_0(Z_1),\Psi_0(L_1))&=& b(Z_1-s_1Z_3, L_1+t_1A +(\sigma(L_1,L_2)+\frac 12 t_1t_2)Z_2)\\
&=& A+t_1Z_1-s_1\sqrt2Z_2+\sqrt2\sigma(L_1,L_2)Z_3
\end{eqnarray*}
and
$$\Psi_0(b(Z_1,L_1))=\Psi_0(A)=A+t_1Z_1+t_2Z_2+(t_3-s_1t_1-s_2t_2)Z_3.$$
Hence, if (B1) holds, then $b(\Psi_0(Z_1),\Psi_0(L_1))= \Psi_0(b(Z_1,L_1))$ is equivalent to (B2).
Moreover, a direct calculation shows that (B1) and (B2) imply $b(X,Y)= \Psi_0(b(X,Y))$ for all other combinations of basis vectors $X$ and $Y$.
\qed
\begin{pr}\label{Pisom}
Let $(\fl,\theta_{\fl},b_\fm)$ and $(\fl',\theta_{\fl'},b_{\fm'})$ be Lie algebras with $B$-structure. Let $\fa$ and $\fa'$ be pseudo-Euclidean spaces with distinguished time-like unit vectors $A$ and $A'$, respectively.  Furthermore, let $(\alpha,\gamma)$, $(\alpha',\gamma')\in\ZQb$ be admissible. 
Then $\fd:=\fd_{\alpha,\gamma}(\fl,\theta_\fl,b_\fm,\fa)$ and $\fd':=\fd_{\alpha',\gamma'}(\fl',\theta_{\fl'},b_{\fm'},\fa')$ are isomorphic as symmetric triples with $\G$-structure if and only if there are an isomorphism $S:(\fl,\theta_{\fl},b_\fm)\rightarrow(\fl',\theta_{\fl'},b_{\fm'})$ and an isometry $U:\fa'\rightarrow \fa$ satisfying $U(A')=A$ such that $(S,U)^*[\alpha',\gamma']=[\alpha,\gamma]\in \HQb$.
\end{pr}
\begin{re}{\rm
In the above proposition the map $$(S,U)^*:{\cal H}_Q^2(\fl',\theta_{\fl'},b_{\fm'},\fa')\rightarrow \HQb$$ is well defined. Indeed, it is easy to check, that $(S,U)^*$ maps ${\cal Z}_Q^2(\fl',\theta_{\fl'},b_{\fm'},\fa')$ to $\ZQb$. Moreover, assume that  
$(\alpha'_1,\gamma_1')$ and $(\alpha'_2,\gamma'_2)$ are equivalent elements of ${\cal Z}_Q^2(\fl',\theta_{\fl'},b_{\fm'},\fa')$. Let ${\cal N}\subset \Aut(\fl,\theta_\fl)$ and ${\cal N}'\subset \Aut(\fl',\theta_{\fl'})$ be the subgroups defined by (\ref{EcN}). Then $(\alpha_1',\gamma_1')=((S_0')^*\alpha'_2,(S_0')^*\gamma'_2)(\tau,\sigma)$ for some $S_0'\in{\cal N}'$ and $(\tau,\sigma)\in{\cal C}_Q^1(\fl',\fa')_+$ that satisfy Properties (B1) and (B2). If we put $S_0:=S^{-1} S_0' S$, then $S_0$ is in ${\cal N}$ and 
$$(S,U)^*(\alpha_1',\gamma_1')=(S_0^*(S,U)^*(\alpha_2',\gamma_2'))\cdot (S,U)^*(\tau,\sigma).$$ 
Hence it suffices to show that $S_0\in{\cal N}$ and $(S,U)^*(\tau,\sigma)\in \CQtheta$ satisfy (B1) and (B2). As for (B1), one checks first that  $\tr (\proj_\fm S_0(b_\fm(L,\cdot)))=\tr (\proj_{\fm'} S'_0(b_{\fm'}(SL,\cdot)))$, then (B1) follows easily. Condition (B2) is easy to check.
}\end{re}
{\sl Proof of Prop.~\ref{Pisom}.} Let us first assume that there exist an isomorphism $S:(\fl,\theta_{\fl},b_\fm)\rightarrow(\fl',\theta_{\fl'},b_{\fm'})$ and an isometry $U:\fa'\rightarrow \fa$ satisfying $U(A')=A$ such that $(S,U)^*[\alpha',\gamma']=[\alpha,\gamma]$. Then $\fd_{(S,U)^*(\alpha',\gamma')}(\fl,\theta_\fl,\fb_\fm,\fa)$ and $\fd$ 
are equivalent, thus isomorphic. Furthermore, it can be easily checked that  $$\Psi_1:\fd_{(S,U)^*(\alpha',\gamma')}(\fl,\theta_\fl,\fb_\fm,\fa)\rightarrow \fd',\quad \Psi_1(Z+A+L)=(S^*)^{-1}(Z)+U^{-1}(A)+S(L)$$
for $Z\in\fl^*,\ A\in\fa,\ L\in\fl$ is an isomorphism.

Now let $\Psi:\fd\rightarrow \fd'$ be an isomorphism. In particular, $\Psi$ is an isomorphism of symmetric triples. Since $(\alpha,\gamma)$ and $(\alpha',\gamma')$ are admissible, we can apply Proposition 6.1. in \cite{KO2}. Hence there is an isomorphism of triples $(\tilde S,U)(\fl,\theta_\fl,\fa)\rightarrow (\fl',\theta_{\fl'},\fa')$ such that $(\tilde S,U)^*(\alpha',\gamma')=(\alpha,\gamma)(\sigma,\tau)$ for some $(\sigma,\tau)\in \CQtheta$. The maps $\tilde S:\fl\rightarrow \fl'$ and  $U^{-1} : \fa \rightarrow \fa'$ are induced by $\Psi$, which
maps the canonical isotropic ideal $\fri=\fl^*$ of $\fd$ to the canonical isotropic ideal $\fri'=(\fl')^*$ and is therefore compatible with the filtrations $\fl^*\subset \fl^*\oplus\fa\subset \fd$ and  $(\fl')^*\subset (\fl')^*\oplus\fa'\subset \fd'$. Furthermore, $\Psi|_{\fl^*}=(\tilde S^*)^{-1}:\fl^*\rightarrow (\fl')^*$.

As above, $\fl=\fm\oplus\fn$, $\fl'=\fm'\oplus\fn'$. Then $(\tilde S^*)^{-1}(\fn^*)=(\fn')^*$ since $\Psi$ must map $b(\fri,\fri)$ to $b'(\fri',\fri')$. Consequently, $\tilde S(\fm)=\fm'$.  Moreover, 
\begin{eqnarray}\label{EtS}
\proj_{\fn'}\tilde S(b_\fm(L,L'))&=&\proj_{\fn'}\Psi (b(L,L'))= \proj_{\fn'}(b'(\Psi (L),\Psi(L'))) \nonumber\\
&=&b'(\proj_{\fm'} \Psi(L),\proj_{\fm'} \Psi(L'))=b_{\fm'}(\tilde S(L),\tilde S(L'))
\end{eqnarray}
for all $L,L'\in\fm$. 
Define $S:(\fl,\theta_\fl)\rightarrow (\fl',\theta_{\fl'})$ and $S_0:(\fl,\theta_\fl)\rightarrow (\fl,\theta_\fl)$ by $\tilde S=S\circ S_0$ and  
$$S(\fm)=\fm',\ S(\fn)=\fn',\ S_0|_\fm=\Id_\fm,\ \proj_\fn S_0|_\fn=\Id_{\fn},\ S_0|_{\fl_+}=\Id_{\fl_+}.$$ 
Then $S_0$ is an automorphism of $(\fl,\theta_\fl)$ by Proposition~\ref{Pl}. Hence $S$ is also an automorphism of $(\fl,\theta_\fl)$. Thus $S_0\in{\cal N}$ and, by (\ref{EtS}), $S:(\fl,\theta_{\fl},b_\fm)\rightarrow(\fl',\theta_{\fl'},b_{\fm'})$ is an isomorphism.  
Since the $\G$-structures on $\fd$ and $\fd$ define orientations on $\fd_-/\fl^*_-$ and $\fd'_-/(\fl')^*_-$, respectively, we obtain $U(A')=A$. Moreover,
$$(\alpha,\gamma)(\sigma,\tau)=(\tilde S,U)^*(\alpha',\gamma')=(S_0,\Id)^*(S,U)^*(\alpha',\gamma'),$$
which proves the claim.
\qed
\subsection{Computation of $\HQb_0$}
In the following we will consider the subset
$\HQb_0\subset\HQb$ of those elements $[\alpha,\beta]\in\HQb$ for which $(\alpha,\gamma)$ is admissible and indecomposable as an element of $\ZQtheta$. This set is acted upon by the group 
$$G:=\Aut(\fl,\theta,b_\fm)\times \Aut(\fa,\ip_\fa,A),$$
where $\Aut(\fa,\ip_\fa,A)=\{U\in \grO(\fa,\ip_\fa),\ U(A)=A\}$.

According to Propositions \ref{Pl} and \ref{Pisom},
in order to give a classification of symmetric triples with $\G$-structure it remains to compute 
$\HQb_0/G$ for $\fl=\fg_{4,1}$ and $\fl=\RR\oplus \fh(1)$ with $\theta_\fl$ and $b_\fm$ as given in Proposition~\ref{Pl}.

For the moment let us forget about the additional structures $\theta_\fl$ and $b_\fm$ and consider $\fg_{4,1}$ and $\RR\oplus \fh(1)$ just as Lie algebras. In \cite{KO1} we introduced the quadratic cohomology set $\HQ$ for a Lie algebra $\fl$ and an arbitrary orthogonal $\fl$-module $\fa$. In \cite{Knil} we determined this cohomology for $\fl\in\{\fg_{4,1},\, \RR\oplus \fh(1)\}$ and trivial $\fl$-modules $\fa$. Let us recall some intermediate results of this computation, which will be useful in the following. Let ${\cal Z}_Q^2(\fl,\fa)$ be the set of quadratic cocycles of a Lie algebra $\fl$ with coefficients in the pseudo-Euclidean vector space $\fa$ considered as a trivial $\fl$-module. If a cocycle $(\alpha,\gamma)\in \ZQtheta$ is admissible, then it is balanced as an element of ${\cal Z}_Q^2(\fl,\fa)$ in the sense of \cite{Knil}.  Moreover, if $(\alpha,\gamma)\in {\cal Z}_Q^2(\fl,\fa)$ is indecomposable and contained in $\ZQtheta$, then it is indecomposable as an element of $\ZQtheta$.

In the following $Z_1,Z_2,Z_3,Z_B\in\fl^*$ will denote the dual basis to $L_1,L_2,L_3,B\in\fl$.
\begin{pr}\label{Padm}
\begin{enumerate}
\item Let $\fl$ be the Lie algebra $\fg_{4,1}=\{[L_2,L_3]=B, [B,L_3]=L_1\}$.
\begin{enumerate}
\item If $(\alpha,\gamma)\in {\cal Z}_Q^2(\fl,\fa)$ is balanced, then $\alpha(L_1,\cdot)\not=0$.
\item If $\fa\in\{\RR^{1,1},\RR^{2,0}\}$ and $A, A_1$ is an orthonormal basis of $\fa$, then the cocycle $(\alpha_1,s\gamma_1)\in{\cal Z}_Q^2(\fl,\fa)$ with
\begin{equation}\label{Eag}
\alpha_1= (Z_1\wedge Z_3)\otimes A_1 + (Z_2\wedge Z_B)\otimes A,\quad
\gamma_1 = Z_B\wedge Z_1\wedge Z_3 
\end{equation}
is balanced and indecomposable for all $s\in\RR$.
\end{enumerate}
\item Let $\fl$ be the Lie algebra $\RR\oplus \fh(1)=\RR\cdot L_1\oplus \{[L_2,L_3]=B\}$. 

If $\fa\in\{\RR^{1,1},\RR^{2,0}\}$ and $A, A_1$ is an orthonormal basis of $\fa$, then $(\alpha_2,\gamma_2)\in{\cal Z}_Q^2(\fl,\fa)$ with
\begin{equation}\label{Eag2}
\alpha_2= (Z_1\wedge Z_3)\otimes A_1 + (Z_B\wedge Z_3)\otimes A,\quad
\gamma_2 = Z_B\wedge Z_1\wedge Z_2 
\end{equation}
is balanced and indecomposable. If $\fa=\RR^{1,0}$, then $(\alpha_3,\gamma_3)\in{\cal Z}_Q^2(\fl,\fa)$ with
\begin{equation}\label{Eag3}
\alpha_3= (Z_B\wedge Z_3)\otimes A,\quad
\gamma_3 = Z_B\wedge Z_1\wedge Z_2 
\end{equation}
is balanced and indecomposable.
\end{enumerate}
\end{pr}
\proof \cite{Knil}, Prop.~4.4 for $\fg_{4,1}$ and Prop.~4.5 for $\fh(1)\oplus\RR$. \qed

\begin{pr} If $(\fl,\theta_\fl,b_\fm)$ is defined as in item 1 of Prop.~\ref{Pl}, then the orbit space $\HQb_0/G$ is not empty if and only if $\fa$ is isomorphic to $\RR^{1,1}$ or $\RR^{2,0}$. For $\fa\in\{\RR^{1,1},\RR^{2,0}\}$ the elements of $\HQb_0/G$ are represented by 
$[\alpha,\gamma_{t}]$, $t\in\RR$, with 
\begin{eqnarray*}
\alpha&=&(-\sqrt 2 Z_B\wedge Z_2 )\otimes A + (Z_1\wedge Z_3)\otimes A_1, \\
\gamma_{t}&=&tZ_B\wedge Z_1\wedge Z_3,
\end{eqnarray*}
where $A_1$ is a fixed unit vector in $\fa_+$. 
\end{pr}

\proof Let $[\alpha,\gamma]$ be in $\HQb_0$. Then $\alpha(L_1,L_2)=0$. Condition (Z2) gives 
$$[B,\sqrt2 L_3]=\langle A,\alpha(B,L_2)\rangle 
L_1-\langle A,\alpha(B,L_1)\rangle L_2=\sqrt2 L_1,$$
hence $\alpha(B,L_2)=-\sqrt 2 A$ and $\alpha(B,L_1)=0$. Because of admissibility $\alpha(L_1,L_3)$ cannot vanish, see Prop.~\ref{Padm}, 1.~(a). Moreover, $\alpha(L_1,L_3)\not=0$ must span $\fa_+$ because of indecomposability. 
Let $S_0\in{\cal N}$ be given by $S_0(L_3)=s_1L_1+s_2L_2+L_3$. Then
\begin{equation}\label{ES0}
(S_0^*\alpha+d\tau)=\alpha- (Z_2\wedge Z_3)\otimes \tau(B)-2\sqrt2 s_2 (Z_B\wedge Z_3)\otimes A 
\end{equation}
since $\tau(L_1)=\sqrt 2 s_2A$ by (B1).
Hence we may assume $\alpha(L_2,L_3)=\alpha(B,L_3)=0$, thus 
$$\alpha=(-\sqrt 2 Z_B\wedge Z_2)\otimes A + (s Z_1\wedge Z_3)\otimes A_1$$
for some unit vector $A_1\in\fa_+$ and some $s\in\RR$, $s\not=0$. It is easy to check that each such $\alpha$ satisfies $d\alpha=0$. Moreover,
for different choices of $s$ the quadratic cocycles $(\alpha,\gamma)\in\ZQb$ are not equivalent. Obviously, $\la\alpha\wedge\alpha\ra=0$ and, moreover, $d\gamma=0$ for each $\theta_\fl$-invariant $\gamma \in C^3(\fl)$. 

As for $\gamma$, note that 
\begin{equation}\label{Ehilf}
2\gamma(B,L_1,L_2)=-\la A,\alpha(B,\sqrt2 L_3)\ra=0
\end{equation}
because of Condition (Z3).

Let us now check how we can change $\gamma$ without changing $\alpha$ and $[\alpha,\gamma]$. By (\ref{ES0}), $(S_0^*\alpha+d\tau)=\alpha$ holds if and only if $\tau(B)=0$ and $s_2=0$, which is equivalent to $d\tau=0$.
Let $\tau$ and $S_0$ satisfy these conditions. Then (\ref{Ehilf}) implies $S_0^*\gamma=\gamma$. Hence 
\begin{equation}\label{Eae}
S_0^*\gamma +d\sigma +\la (S_0^*\alpha +\textstyle{\frac12} d\tau )\wedge\tau\ra=\gamma+d\sigma+\la \alpha\wedge \tau\ra.
\end{equation}
Since 
$$d\sigma=\sigma(L_1,L_2)\cdot Z_B\wedge Z_2\wedge Z_3=-(1/\sqrt 2)\la \tau(L_3),A\ra\cdot Z_B\wedge Z_2\wedge Z_3$$
by Condition (B2) and 
$$\la \alpha\wedge \tau\ra=\la\alpha(B,L_2),\tau(L_3)\ra\cdot Z_B\wedge Z_2\wedge Z_3=\la -\sqrt2 A,\tau(L_3)\ra\cdot Z_B\wedge Z_2\wedge Z_3$$
the right hand side of (\ref{Eae}) equals
$$\gamma-(3/\sqrt 2)\la \tau(L_3),A\ra\cdot Z_B\wedge Z_2\wedge Z_3.$$
Hence we may assume $\gamma(B,L_2,L_3)=0$. Together with (\ref{Ehilf}) this implies
$$\gamma=tZ_B\wedge Z_1\wedge Z_3$$ for some $t\in\RR$.
It remains to decide which of these $(\alpha,\gamma)$ are admissible and indecomposable and to divide by~$G$.

It holds that $S \in \Aut(\fl,\theta_\fl,b_\fm)$ and $U \in\Aut(\fa,\ip_\fa,A)$ if and only if 
$$ S(L_1)=a^3 L_1,\quad S(L_2)=b L_1+(1/a)\cdot L_2,\quad S(L_3)=a^2 L_3,\quad S(B)=aB, $$
$$U(A)=A,\quad U(A_1)=\delta A_1,\quad \delta=\pm1,$$
for some $a,b\in\RR$, $a\not=0$. Hence $[\alpha,\gamma]$ is in the same orbit as $[\alpha',\gamma']$ with
\begin{eqnarray*}
\alpha'&=&(-\sqrt 2 Z_B\wedge Z_2 )\otimes A + (Z_1\wedge Z_3)\otimes A_1, \\
\gamma'&=&t'Z_B\wedge Z_1\wedge Z_3
\end{eqnarray*}
for some $t'\in\RR$. For different choices of the parameters $t'$ the equivalence classes $[\alpha',\gamma']$ belong to different $G$-orbits. 
As a cohomology class in ${\cal H}_Q^2(\fl,\fa)$, $[\alpha',\gamma']$ equals 
$$[(\sqrt2 Z_2\wedge Z_B)\otimes A+(Z_1\wedge Z_3)\otimes A_1,t'Z_B\wedge Z_1\wedge Z_3],$$
which is in the same $\Aut(\fl)\times\grO(\fa,\ip_\fa)$-orbit as $[\alpha_1,s\gamma_1]$ for some $s\in\RR$ with $\alpha_1$ and $\gamma_1$ as in (\ref{Eag}). Hence $(\alpha',\gamma')$ admissible and indecomposable. \qed
\begin{pr}If $(\fl,\theta_\fl,b_\fm)$ is defined as in item 2 of Prop.~\ref{Pl}, then the orbit space $\HQb_0/G$ is not empty if and only if $\fa$ is isomorphic to $\RR^{1,1}$, $\RR^{2,0}$ or $\RR^{1,0}$. Put 
\begin{equation}\gamma_0=(1/\sqrt2)\cdot Z_B\wedge Z_1\wedge Z_2.\label{Egrs}
\end{equation}
For $\fa\in\{\RR^{1,1},\RR^{2,0}\}$ we have $\HQb_0/G=\{[\bar \alpha,\gamma_0]\}$  with 
\begin{equation}
\bar\alpha=( Z_B\wedge Z_3 )\otimes A + (Z_1\wedge Z_3)\otimes A_1, \label{Eak}
\end{equation}
where $A_1$ is a fixed unit vector in $\fa_+$. 

For $\fa=\RR^{1,0}$ we have $\HQb_0/G=\{[\alpha_0,\gamma_0]\}$,   where 
$$
\alpha_0=( Z_B\wedge Z_3 )\otimes A.
$$
\end{pr}
\proof Let $[\alpha,\gamma]$ be in $\HQb_0$. Then $\alpha(L_1,L_2)=0$.  Condition (${\rm Z2}$) gives 
$$[B,\sqrt2 L_3]=\langle A,\alpha(B,L_2)\rangle 
L_1-\langle A,\alpha(B,L_1)\rangle L_2=0,$$
hence $\alpha(B,L_1)=\alpha(B,L_2)=0$.  
Let $S_0\in{\cal N}$ be given by $S_0(L_3)=s_1L_1+s_2L_2+L_3$. Then
\begin{equation}\label{Es0z}
(S_0^*\alpha+d\tau)=\alpha - (Z_2\wedge Z_3)\otimes \tau(B).
\end{equation}
Thus we may assume $\alpha(L_2,L_3)=0$. Moreover, $\alpha(L_1,L_3)$ must span $\fa_+$. Thus we have 
$$\alpha=(t Z_B\wedge Z_3)\otimes A + (t' Z_1\wedge Z_3)\otimes A_1$$
for some $t'\not=0$ and a unit vector $A_1\in\fa_+$, if $\fa_+\not=0$ and
$$\alpha=(t Z_B\wedge Z_3)\otimes A $$
if $\fa_+=0$. In both cases $t\not=0$ because of indecomposability. Obviously, each such $\alpha$ satisfies $d\alpha=0$ and $\la\alpha\wedge\alpha\ra=0$. For different choices of $t$ and $t'$ the quadratic cocycles are in different equivalence classes. Condition (Z3) gives
$$2\gamma(B,L_1,L_2)=-\la A,\alpha(B,\sqrt2 L_3)\ra=-\la A,\sqrt2 t A\ra=\sqrt2 t.$$ 
This is the only condition for $\gamma$ since $d\gamma=0$ for each $\theta_\fl$-invariant $\gamma \in C^3(\fl)$. Let us now check how we can change $\gamma$ without changing $\alpha$ and $[\alpha,\gamma]$. By (\ref{Es0z}) $(S_0^*\alpha+d\tau)=\alpha$ holds if and only if $\tau(B)=0$. Moreover, (B1) implies $\tau(L_1)=\sqrt2 s_2$ and $\tau(L_2)=-\sqrt2 s_1$. Let $\tau$ and $S_0$ satisfy these conditions. Then 
\begin{equation}\label{Egs}
S_0^*\gamma +d\sigma +\la (S_0^*\alpha +\textstyle{\frac12} d\tau )\wedge\tau\ra=S_0^*\gamma+\la \alpha\wedge \tau\ra.
\end{equation}
Since
$$ S_0^*\gamma=\gamma+(s_2t/\sqrt2)\cdot Z_B\wedge Z_1\wedge Z_3 - (s_1t/\sqrt2)\cdot  Z_B\wedge Z_2\wedge Z_3$$ and
$$\la\alpha\wedge\tau\ra=\sqrt2 s_2tZ_B\wedge Z_1\wedge Z_3-\sqrt 2 s_1tZ_B\wedge Z_2\wedge Z_3.$$  
Equation (\ref{Egs}) implies that we may assume $\gamma=(t/\sqrt2)\cdot Z_B\wedge Z_1\wedge Z_2$.

It holds that $S \in \Aut(\fl,\theta_\fl,b_\fm)$ if and only if 
$$ S(L_1)=a L_1,\quad S(L_2)=b L_2+xL_1,\quad S(L_3)=ab L_3,\quad S(B)=ab^2B $$
and $U \in\Aut(\fa,\ip_\fa,\theta_\fa,A)$ if and only $U=\Id$ in case $\fa=\RR^{1,0}$ and 
$$U(A)=A,\quad U(A_1)=\delta A_1,\quad \delta=\pm1$$
if $\fa\in\{\RR^{2,0},\RR^{1,1}\}$. Hence, if $\fa\in\{\RR^{2,0},\RR^{1,1}\}$, then $[\alpha,\gamma]$ is in the same orbit as $[\bar\alpha,\gamma_0]$ with $\bar \alpha$ and $\gamma_0$  as defined in (\ref{Eak}) and (\ref{Egrs}).  As a cohomology class in ${\cal H}_Q^2(\fl,\fa)$, $[\bar \alpha,\gamma_0]$ is in the same $\Aut(\fl)\times\grO(\fa,\ip_\fa)$-orbit as $[\alpha_{2},\gamma_{2}]$  with $\alpha_2$ and $\gamma_2$ as in (\ref{Eag2}). Hence $(\bar \alpha,\gamma_0)$ is admissible and indecomposable. 

If $\fa=\RR^{1,0}$, then $[\alpha,\gamma]$ is in the same $G$-orbit as $[\alpha_0,\gamma_{0}]$.
As a cohomology class in ${\cal H}_Q^2(\fl,\fa)$, $[\alpha_0,\gamma_0]$ is in the same $\Aut(\fl)\times\grO(\fa,\ip_\fa)$-orbit as $[\alpha_3,\gamma_3]$ for $\alpha_3,\gamma_3$ as in (\ref{Eag3}). Hence $(\alpha_0,\gamma_0)$ is admissible and indecomposable. \qed

\subsection{Classification result}\label{Sfinal}
Recall that for a given Lie algebra $\fl$ with involution $\theta_\fl$, a pseudo-Euclidean space $(\fa,\ip_\fa)$ with isometric involution $\theta_\fa$, a suitable 2-cocycle $\alpha\in Z^2(\fl,\fa)$ and a suitable 3-form $\gamma\in\bigwedge^3\fl^*$, we can define a symmetric triple $\dd=(\fd,\theta,\ip)$ by $\fd=\fl^*\oplus\fa\oplus\fl$ (as a vector space), $\theta=\theta_\fl^*\oplus\theta_\fa\oplus\theta_\fl$, 
$[\fl^*\oplus\fa,\fl^*\oplus\fa] =0$,
\begin{eqnarray*}
\ [L,L'] &=& \gamma(L,L',\cdot) +\alpha(L,L')+[L,L']_\fl\\
\ [L,A+Z] &=&  - \langle A,\alpha(L,\cdot)\rangle+ \ad ^*(L)(Z)
\end{eqnarray*}
and 
$$
 \langle Z+A+L,Z'+A'+L'\rangle= \langle A,A'\rangle_\fa
+Z(L') +Z'(L)$$
for $Z,\,Z'\in \fl^*$, $A,\,A'\in \fa$ and
$L,\,L'\in \fl$.

In the following theorem $Z_1,Z_2,Z_3,Z_B\in\fl^*$ is the dual basis to $L_1,L_2,L_3,B\in\fl$.
\begin{theo} If $(\fg,\theta,\omega,\ip)$ is an indecomposable symmetric triple with $\G$-struc\-ture, then it is isomorphic to  exactly one of the  symmetric triples with $\G$-structure $(\dd, \omega)$ for the following data $\fl,\fa,\alpha,\gamma$ and $\omega$:
\begin{enumerate}
\item $\fl=\fg_{4,1}=\{B=[L_2,L_3]_\fl,\ [B,L_3]_\fl=L_1\}$, 
$\fl_-=\Span\{L_1,L_2,L_3\}$,  $\fl_+=\RR\cdot B$,\\
$\fa\in\{\RR^{1,1},\RR^{2,0}\}$ with fixed orthonormal basis $A, A_1$,\\
$\alpha=(-\sqrt 2 Z_B\wedge Z_2 )\otimes A + (Z_1\wedge Z_3)\otimes A_1,$ \\
$\gamma=tZ_B\wedge Z_1\wedge Z_3$, $t\in\RR$; or
\item $\fl=\RR\oplus\fh(1)=\RR\cdot L_1\oplus\{B=[L_2,L_3]_\fl\}$,
$\fl_-=\Span\{L_1,L_2,L_3\}$,  $\fl_+=\RR\cdot B$,
\begin{enumerate}
 \item $\fa\in\{\RR^{1,1},\RR^{2,0}\}$ with fixed orthonormal basis $A, A_1$,\\
 $\alpha=( Z_B\wedge Z_3 )\otimes A + (Z_1\wedge Z_3)\otimes A_1$, \\
$\gamma=(1/\sqrt2)\cdot Z_B\wedge Z_1\wedge Z_2$; or
\item $\fa=\RR^{1,0}$,\\
$\alpha=( Z_B\wedge Z_3 )\otimes A$,  $\gamma=(1/\sqrt2)\cdot Z_B\wedge Z_1\wedge Z_2$; 
\end{enumerate}
\end{enumerate}
and
$$\omega=\sqrt2 (\sigma^{127}+\sigma^{356})-\sigma^4\wedge(\sigma^{15}+\sigma^{26}-\sigma^{37}),$$
where  $\sigma^1,\dots,\sigma^7$ is the dual basis to $Z_1,Z_2,Z_3,A,L_1,L_2,L_3$.

All listed symmetric triples $\dd$ are indecomposable and pairwise non-isomorphic. 
 \end{theo}
 
 A direct consequence is the following corollary, which can also be deduced already from Prop.~\ref{Pl}.

\begin{co} If an indecomposable symmetric space $(M,g)$ of signature $(4,3)$ admits a parallel $\G$-structure, then its transvection group is nilpotent, and its holonomy group is three-dimensional and abelian.
\end{co}

\end{document}